%
%
%

\input amstex
\documentstyle{amsppt}
\PSAMSFonts

\pageheight{18.5cm}
\magnification=\magstep1
\frenchspacing


\loadbold

\font\smallbf=cmb10 at 8pt
\font\small=cmr10 at 8pt

\def\phi{\varphi}

\def\nin{\newline\indent}
\def\name#1{{\smc #1\/}}
\def\CC{{\Bbb C}}
\def\FF{{\Bbb F}}

\def\PP{{\Bbb P}}

\def\ZZ{{\Bbb Z}}

\def\CCC{{\Cal C}}

\def\OOO{{\Cal O}}

\def\se#1,#2,#3;{(ab)^{#1}\, a_x^{\,#2}b_x^{\,#3}}

\def\see#1,#2,#3;{(ab)^{#3} (bc)^{#1} (ca)^{#2}\,
\ell_a^{p-#2-#3}\ell_b^{q-#1-#3}\ell_c^{r-#1-#2}}

\def\seep#1,#2,#3;{(ab)^{#3} (bc)^{#1} (ca)^{#2}\,
\ell_a^{p'-#2-#3}\ell_b^{q'-#1-#3}\ell_c^{r'-#1-#2}}

\def\seed#1,#2,#3;{(ab)^{#3} (bc)^{#1} (ca)^{#2}\,
a_x^{p_a-#2-#3}b_x^{p_b-#1-#3}c_x^{p_c-#1-#2}}

\def\nse#1,#2,#3,#4,#5,#6;{(#1 #2)(#3 #4){#5}_x{#6}_x}
\def\nsee#1,#2,#3,#4,#5,#6;{(#1 #2)(#3 #4)#5_x#6_x}
\def\nsse#1,#2,#3,#4,#5,#6,#7,#8;{(#1 #2)(#3 #4)(#5 #6)(#7 #8)}

\def\tv#1,#2,#3;{[#1,#2]_{#3}}

\TagsOnRight

\input xy
\xyoption{all}

\topmatter
\title Geometry of orbit closures for the representations associated to gradings of Lie algebras of types $E_7$.
\endtitle
\date December 29, 2011 \enddate
\author Witold Kra\'skiewicz and Jerzy Weyman \endauthor
\address Nicholas Copernicus University
\nin Toru\'n , Poland
\endaddress
\email wkras\@mat.uni.torun.pl \endemail
\address Department of Mathematics, Northeastern University
\nin 360 Huntington Avenue,  BOSTON,  MA 02115, USA \endaddress
\email j.weyman\@neu.edu \endemail

\thanks The  second author was partially supported by NSF grant DMS-0600229
\endthanks
\rightheadtext{Geometry of orbit closures for $E_7$}
\abstract
This paper is a continuation of \cite{KW11a}.  We  investigate the orbit closures for the class of representations of simple algebraic groups
associated to various gradings on the simple Lie algebra of type $E_7$. The methods for classifying the
orbits for these actions were developed by Vinberg \cite{V75}, \cite{V87}.
We give the orbit descriptions, the degeneration partial orders, and indicate normality of the orbit closures.
We also investigate the rational singularities, Cohen-Macaulay and Gorenstein properties for the orbit closures.
We give the information on the defining ideals of orbit closures.
The corresponding results for the Lie algebras of types $E_6$, $F_4$, $G_2$ were given in \cite{KW11a}.

\endabstract
\endtopmatter

\document
\def\iz#1:#2.{\noindent\S\ignorespaces#1\dotfill#2\par}

\head  Introduction \endhead

The irreducible representations of semi-simple groups with finitely many orbits were classified by Kac in \cite{K82} (with some corrections in \cite{DK85}).
They correspond (with very few exceptions) to certain gradings on the root systems, and to the corresponding $\theta$ groups.
In Kac's paper the list of these representations appears in Table II (there are the other tables III, IV classifying so-called visible representations).
We refer to these representations as representations of type I.

 The representations of type I are parametrized by the pairs $(X_n, \alpha_k)$ where $X_n$ is a  Dynkin diagram with a distinguished node
 $x\in X_n$. This data defines a grading
 $${\goth g}=\oplus_{i=-s}^s {\goth g}_i$$
 of a simple algebra $\goth g$ of type $X_n$ such that the Cartan subalgebra $\goth h$ is contained in ${\goth g}_0$ and the root space ${\goth g}_\beta$ is contained in ${\goth g}_i$ where $i$ is the coefficient of the simple root $\alpha$ corresponding to the node $x$  in the expression for $\beta$ as a linear combination of simple roots. The representation corresponding to $(X_n ,x)$ is the ${\goth g}_1$ with the action of the group $G_0\times \CC^*$ where $G_0$ is the adjoint group corresponding to ${\goth g}_0$ and $\CC^*$ is the copy of $\CC^*$ that occurs in maximal torus of $G$ (the adjoint group corresponding to $\goth g$)  but not in maximal torus of $G_0$.

The orbit closures for the representations of type I were described in two ways by Vinberg in \cite{V75}, \cite{V87}.
The first description states that the orbits are the irreducible components of the intersections of the nilpotent orbits in $\goth g$
with the graded piece ${\goth g}_1$.
In the second paper Vinberg gave a more precise description in terms of the support subalgebras which are graded Lie subalgebras
of the graded Lie algebra $\goth g =\oplus {\goth g}_i$

In this paper we  concentrate on the cases when $X_n$ is equal to $E_7$ .
The corresponding results for the Lie algebras of types $E_6$, $F_4$, $G_2$ were given in \cite{KW11a}.

The main result of the paper is the calculation of the Hilbert polynomials of the normalizations of the coordinate rings
of the orbit closures. We are also able to decide the normality, Cohen-Macaulay and Gorenstein
properties of the orbit closures as well as rational singularities property.

In some cases we are also able to describe explicitly the free resolutions of the coordinate rings of orbit closures as modules over the
coordinate ring of the representation itself.  We list the terms of these resolutions in the corresponding sections.
These calculations are being carried out by Federico Galetto and will be published elsewhere.

The idea of such calculations is as follows.  First, one gets the suggested terms of the complex from the geometric method.
Then one tries to use the interactive  Macaulay 2 calculations using Buchsbaum-Eisenbud exactness criterion
 \cite {BE73}, see the final remark of section 3 of \cite{KW11a}.

The paper is organized as follows. In section 1 we introduce the necessary notation.
In the remaining sections we present the data for the Lie algebra of type $E_7$.

The orbits were calculated first by hand but then the calculations were checked using the program \cite{dG11} kindly provided by Willem de Graaf.
The dimensions of the orbits were calculated using computer routines written by Jason Ribeiro. Jason also wrote a very useful
 python package \cite{R10} for calculating  Euler characteristics of the bundles involved.

 The bulk of the calculations was done using  more complicated roots and weight programs (written by the first author) which searched through all parabolic subgroup submodules
 for the representations in question. Then the program calculated needed Euler characteristics and Hilbert polynomials.

 The data are organized as follows. For each representation we start with several  tables.
 First there is a general table with the number of the orbit, the type of the support algebra $\goth s$ and the dimension of the orbit.
 This is followed by  tables with the geometric description of the orbits.
 The numerical data table includes the degree and the numerator of the Hilbert series of the coordinate ring of the normalization of the orbit closure. The denominator is always $(1-t)^{codim}$.
The final table indicates the singularities data i.e. the information on whether the orbit closure and its normalization is spherical, normal, Cohen-Macaulay, has rational singularities and is Gorenstein.

In the final section we list some general conclusions about the orbit closures in the representations we deal with.

\bigskip\noindent
{\smallbf Acknowledgment:} {\small
Both authors thank Federico Galetto, Willem de Graaf and Jason Ribeiro
for very useful computer programs.
The second author would like to thank Joe Landsberg and Steven V. Sam
 for interesting conversations and pointing out various errors in earlier versions of this paper. }

\bigskip\bigskip

\head  \S 1. Preliminaries and notation \endhead

Let $X_n$ be a Dynkin diagram and let $\goth g$ be the corresponding
 simple Lie algebra. Let us distinguish a node $x\in X_n$. Let $\alpha_k$ be
a corresponding simple root in the root system $\Phi$ corresponding
to $X_n$. The choice of $\alpha_k$ determines a $\ZZ$- grading on $\Phi$
by letting the degree of a root $\beta$ be equal to the coefficient of
$\alpha_k$ when we write $\beta$ as a linear combination of simple roots. On the level
of Lie algebras this corresponds to a $\ZZ$-grading
$${\goth g}=\oplus_{i\in \ZZ}\ {\goth g}_i .$$

We define the group $G_0 := (G ,G )\times {\CC}^*$ where $(G,G)$ is a
connected semisimple group with the Dynkin diagram $X_n\setminus x$.
A representation of type I is the representation of $G_0$ on ${\goth g}_1$.

We will denote the representation ${\goth g}_1$ by $X_n,\alpha_k)$.

Denoting by ${\goth l}$ the Levi factor ${\goth g}_0$ we have
$${\goth l}={\goth l}'\oplus {\goth z}(\goth l )$$
where ${\goth l}'$ denotes the Lie algebra associated to $X_n$ with the omitted node $x$, and ${\goth z}({\goth l})$ is a one dimensional center of $\goth l$.

In this paper the Lie algebra $\goth g$ is a simple Lie algebra of type $E_7$.

 \name{Vinberg} in \cite{V75}, \cite{V87}  gave two descriptions of the $G_0$-orbits in the representations of type $II$  in terms of
conjugacy classes of nilpotent elements in $\goth g$. We refer to \cite{KW11a}, sections 1 and 2 for the precise statements.

All the orbit closures in the representations we consider have a desingularization by a total space of homogeneous vector bundle
over the appropriate homogeneous space $G/P$. The description of the results we use in this context is given in \cite{KW11a} section 3.

Each orbit closure $\overline{\OOO_x}$ has a desingularization $Z(x,h,y)$ associated to the ${\goth sl}_2$-triple $(x,h,y)$ with $x\in {\goth g}_1$, $h\in{\goth g}_0$, $y\in{\goth g}_{-1}$ described in the section 4 of \cite{KW11a}.

In the case the orbit closure $\overline{\OOO}_x$ is not normal, it satisfies the condition of the Remark 4.3 from \cite{KW11a}.

Some of the orbit closures are {\it degenerate} i.e. their singularities come from another orbit closure from a smaller representation of type $(X_n,\alpha_k)$
where $X_n$ is some proper subdiagram of $E_7$. In such situation we can deduce several properties of the bigger orbit closure from the properties of the smaller one.
We collected necessary facts in \cite{KW11a} section 5.

\proclaim{Remark 1.1. (proving normality and rational singularities)}
In the cases we consider below we claim the general fact that the normalizations of orbit closures have rational singularities and we list the normal orbit closures.This is done as in \cite{KW11a}. The analogues of Proposition 3.7. and Proposition 3.8 from \cite{KW11a} are true for all orbits for representations related to the gradings of Lie algebra  $E_7$.
\endproclaim

\bigskip\bigskip

\head \S 2. The case $(E_7, \alpha_1)$. \endhead

 The representation in question is $X = V(\omega_5 , D_6 )$, a half-spinor
representation for the group
${ G_0}=Spin (12)\times\CC^*$. Here $\omega_5 = ({1\over 2}, {1\over 2}, {1\over 2}, {1\over 2}, {1\over 2}, {1\over 2})$.
$dim(X)=32$.
The weights of $X$ with respect to $Spin(12)$ are vectors in 6 dimensional space, with coordinates equal to $\pm{1\over 2}$, with even number of negative coordinates.

The graded Lie algebra of type $E_7$ is
$${\goth g}(E_7)= {\goth g}_{-2}\oplus {\goth g}_{-1}\oplus {\goth g}_0\oplus {\goth g}_1\oplus {\goth g}_2$$
with ${\goth g}_0= \CC\oplus {\goth so}(12)$,
${\goth g}_1=V(\omega_5 ,D_6)$, ${\goth g}_2=\CC$.

We denote the weight vectors by the subsets $[I]$ where $I$ is the subset of the set $\lbrace 1,2,3,4,5,6\rbrace$ of even cardinality where the
component of a given weight vector is negative.

The invariant scalar product on ${\goth h}$ restricted to the roots from ${\goth g}_1$ is

$$([I], [J] )=2-{1\over 2} \#((I\setminus J)\cup (J\setminus I)).$$
Notice that the possible scalar products are $2, 1, 0, -1$.

This is another member of ``subexceptional
series" of Landsberg and Manivel \cite{LM01}. The ring of invariants is generated by a
discriminant $\Delta$ of degree 4. There are five  orbits with linear containment diagram


$$
\xy
(15,0)*+{{\Cal O}_{0}}="o0";%
(15,8)*+{{\Cal O}_{1}}="o1";%
(15,16)*+{{\Cal O}_2}="o2";%
(15,24)*+{{\Cal O}_3}="o3";%
(15,32)*+{{\Cal O}_4}="o4";%
(-15,0)*{0};%
(-15,8)*{16};%
(-15,16)*{25};%
(-15,24)*{31};%
(-15,32)*{32};%
{\ar@{-} "o0"; "o1"};%
{\ar@{-} "o1"; "o2"};%
{\ar@{-} "o3"; "o2"};%
{\ar@{-} "o4"; "o3"};%
\endxy
$$

$$\matrix number&{\goth s}&dim&representative\\
0&0&0&0\\
1&A_1&16&[\emptyset ]\\
2&2A_1&25&[\emptyset]+[1234] \\
3&3A_1&31&[\emptyset ]+[1234]+[1256]\\
4&A_2&32&[\emptyset ]+[123456]\
 \endmatrix$$

$$\matrix number&proj.\ picture&tensor\ picture\\
0&0&0\\
1&h.w.\ vector&pure\ spinors&\\
2&&sing(hyperdisc.)\\
3&\tau ({\overline\OOO}_1)&hyperdisc.\\
4&\sigma_2 ({\overline\OOO}_1)&generic \endmatrix$$

The numerical data is as follows:

$$\matrix number&degree&numerator\\
0&1&1\\
1&286&1+16t+70t^2+112t^3+70t^4+16t^5+t^6\\
2&176&1+7t+28t^2+52t^3+52t^4+28t^5+7t^6+t^7\\
3&4&1+t+t^2+t^3\\
4&1&1\\
\endmatrix$$

The singularities data is as follows.

$$\matrix number&spherical&normal&C-M&R.S.&Gor\\
0&yes&yes&yes&yes&yes\\
1&yes&yes&yes&yes&yes\\
2&yes&yes&yes&yes&yes\\
3&yes&yes&yes&yes&yes\\
4&no&yes&yes&yes&yes \endmatrix$$

$\spadesuit$ The orbit closure $\overline{\OOO_3}$ is a hypersurface given
by the invariant of degree 4.

$\spadesuit$ Let us look at the orbit closure $\overline{\OOO_2}$.
The resolution of the coordinate ring is
$$0\rightarrow A(-14)\rightarrow V_{\omega_6}\otimes A(-11)\rightarrow V_{\omega_2}\otimes A(-10)\rightarrow V_{2\omega_1}\otimes A(-8)\rightarrow$$
$$\rightarrow V_{2\omega_1}\otimes A(-6)\rightarrow V_{\omega_2}\otimes A(-4)\rightarrow V_{\omega_6}\otimes A(-3)\rightarrow A$$

\bigskip\bigskip

\head \S 3. The case $(E_7, \alpha_2)$. \endhead

  $X= \bigwedge^3 F$, $F=\CC^7$, ${ G}=GL(F)$. The orbits in this case were
calculated for the first time in the book \cite{Gu64} of Gurevich.

The graded Lie algebra of type $E_7$ is
$${\goth g}(E_7)= {\goth g}_{-2}\oplus {\goth g}_{-1}\oplus {\goth g}_0\oplus {\goth g}_1\oplus {\goth g}_2$$
with $G_0= GL(7)$, ${\goth g}_0= \CC\oplus {\goth sl}(7)$,
${\goth g}_1=\bigwedge^3\CC^7$, ${\goth g}_2=\bigwedge^6\CC^7$.

The weights of ${\goth g}_1$ are $\epsilon_i +\epsilon_j+\epsilon_k$ for $1\le i<j<k\le 7$.
We label thie corresponding weight vector by $[I]$ where $I$ is a cardinality 3 subset of $\lbrace 1,2,3,4,5,6,7\rbrace$.

The invariant scalar product on ${\goth h}$ restricted to the roots from ${\goth g}_1$ is
$$([I], [J])= \delta-1$$
where $\delta =\# (I\cap J)$.

This representation has ten orbits. The containment diagram is
$$
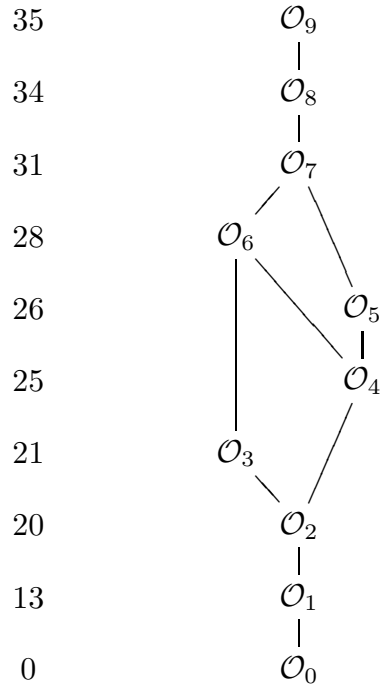
\xy
(15,0)*+{{\Cal O}_{0}}="o0";%
(15,8)*+{{\Cal O}_{1}}="o1";%
(15,16)*+{{\Cal O}_2}="o2";%
(8,24)*+{{\Cal O}_3}="o3";%
(22,32)*+{{\Cal O}_4}="o4";%
(22,40)*+{{\Cal O}_5}="o5";%
(8,48)*+{{\Cal O}_6}="o6";%
(15,56)*+{{\Cal O}_7}="o7";%
(15,64)*+{{\Cal O}_8}="o8";%
(15,72)*+{{\Cal O}_9}="o9";%
(-15,0)*{0};%
(-15,8)*{13};%
(-15,16)*{20};%
(-15,24)*{21};%
(-15,32)*{25};%
(-15,40)*{26};%
(-15,48)*{28};%
(-15,56)*{31};%
(-15,64)*{34};%
(-15,72)*{35};%
{\ar@{-} "o0"; "o1"};%
{\ar@{-} "o1"; "o2"};%
{\ar@{-} "o3"; "o2"};%
{\ar@{-} "o4"; "o2"};%
{\ar@{-} "o4"; "o5"};%
{\ar@{-} "o3"; "o6"};%
{\ar@{-} "o4"; "o6"};%
{\ar@{-} "o5"; "o7"};%
{\ar@{-} "o6"; "o7"};%
{\ar@{-} "o7"; "o8"};%
{\ar@{-} "o8"; "o9"};%
\endxy
$$

 $$\matrix number&\goth s &dim&representative\\
0&0&0&0\\
1&A_1&13&[123]\\
2&2A_1&20&[123]+[145]\\
3&3A_1&21&[123]+[145]+[167]\\
4&3A_1&25&[123]+[145]+[246]\\
5&A_2&26&[123]+[456]\\
6&4A_1&28&[123]+[145]+[167]+[357]\\
7&A_2+A_1&31&[123]+[456]+[147]\\
8&A_2+2A_1&34&[123]+[456]+[147]+[257] \\
9&A_2+3A_1&35&[123]+[456]+[147]+[257]+[367]
  \endmatrix$$

  $$\matrix numberproj.\ picture&tensor\ picture\\
0&0&0\\
1&C(Grass(3,7))&h.w.\ vector\\
2&&tensors\ of\ rank\le 5\\
3&&1-decomposable\ tensors\\
4&\tau(\overline\OOO_1)&$F$-degenerate, hyperdisc.\\
5&\sigma_2 (\overline\OOO_1)&tensors\ of\ rank\le 6\\
6&J(\overline\OOO_1 ,\tau(\overline\OOO_1))&polarizations\ of\ hyperdisc.\ for\ \bigwedge^3 (\CCC^6 )\\
7&\sigma_3 (\overline\OOO_1)&sing(hyperdisc.)\\
8&&hyperdisc. \\
9&&generic  \endmatrix$$

The numerical data is as follows.

$$\matrix number&degree&numerator\\
0&1&1\\
1&462&1+22t+113t^2+190t^3+113t^4+22t^5+t^6\\
2&2394&1+15t+120t^2+428t^3+750t^4+687t^5+316t^6+\\
&&+70t^7+7t^8\\
3&1366&1+14t+105t^2+ 336t^3+ 490t^4+ 336t^5+84t^6\\
4&1792&1+10t+55t^2+192t^3+407t^4+511t^5+385t^6+\\
&&+175t^7+49t^8+7t^9\\
5&735&1+9t+ 45t^2+ 137t^3+ 243t^4+ 216t^5 +84t^6\\
6&1024&1+7t+ 28t^2+ 84t^3+182t^4+ 266t^5+ 252t^6+\\
&&+ 148t^7+ 49t^8+ 7t^9\\
7&210&1+4t+10t^2+20t^3+35t^4+56t^5+49t^6+\\
&&+28t^7+7t^8\\
8&7&1+t+t^2+t^3+t^4+t^5+t^6\\
9&1&1
\endmatrix$$

The singularities data are as follows.

 $$\matrix number&spherical&normal&C-M&R.S.&Gor\\
0&yes&yes&yes&yes&yes\\
1&yes&yes&yes&yes&no\\
2&yes&yes&yes&yes&no\\
3&yes&yes&yes&yes&no\\
4&yes&yes&yes&yes&no\\
5&no&yes&yes&yes&no\\
6&yes&yes&yes&yes&no\\
7&no&yes&yes&yes&no\\
8&no&yes&yes&yes&yes \\
9&no&yes&yes&yes&yes  \endmatrix$$

We will describe in detail the non-degenerate orbit closures in $\bigwedge^3\CC^7$.
These are the orbits $\overline{\OOO_9}$, $\overline{\OOO_8}$. $\overline{\OOO_7}$, $\overline{\OOO_6}$ and $\overline{\OOO_3}$.
The first of these is generic so there is not much to say. We also describe the generic degenerate orbit of tensors of rank $\le 6$.

We use the usual notation. $A=Sym(\bigwedge^3F^*)$ and $(a,b,c,d,e,f,g)$ abbreviates for $S_{a,b,c,d,e,f,g}F^*$.

\bigskip\bigskip

$\spadesuit$  The hyperdiscriminant orbit $\OOO_{8}$.

This is the hypersurface given by the tensors with vanishing hyperdiscriminant.
The orbit closure  ${\overline\OOO}_8=Y_{hw}^\vee$  is characterized
(set-theoretically) by the condition $S_{3^7}F^* = 0$, $S_{3^4 ,2^3}F^*\ne 0$.

Its desingularization lives on $Grass(3,F)$.
We denote by $\Cal R$, $\Cal Q$ the tautological subbundle and factorbundle respectively.
The bundle $\xi$ is
$$\xi = \bigwedge^3 {\Cal R}+{\Cal Q}\otimes\bigwedge^2{\Cal R}.$$
We have $dim\ Z(8)= 4+3\times 6+4\times 3=34$ and as always (see \cite{KW11a}, section 5) $Z(8)$ gives a desingularization of the hyperdiscriminant hypersurface.

The complex $\FF(8)_{\bullet}$ is
$$0\rightarrow (3^7)\rightarrow  (0^7).$$
The hyperdiscriminant $\Delta$ has degree $7$ and it defines a normal hypersurface with rational singularities.

\bigskip\bigskip

$\spadesuit$  The codimension $4$ orbit $\OOO_{7}$.

This orbit closure is the singular locus of the hyperdiscriminant orbit $\overline{\OOO_8}$.

The minimal elements in the bundle $\eta$ describing the desingularization $Z(7)$  are the weights $[1,2,7]$ and $[2,5,6]$.
The bundle $\eta$ is defined over the flag variety $Flag(2, 6; F)$. It has rank $17$, so the dimension of the desingularization is
$17+14=31$ as needed.

One gets a very nice complex describing the resolution of $\CC[N({\overline{\OOO_7}})]$.

 The terms of the complex $\FF(7)_\bullet$ are as follows
$$0\rightarrow (6,5^6)\rightarrow (5^2,4^5)\rightarrow (4,3^5,2)\rightarrow(3^4,2^3)\rightarrow (0^7).$$

The orbit closure is normal.

\bigskip\bigskip

$\spadesuit$ The codimension $7$ orbit $\OOO_{6}$.

The minimal elements in the bundle $\eta$ describing the desingularization $Z(6)$  are the weights $[1,4,7]$ and $[2,3,4]$.
The bundle $\eta$ is defined over the flag variety $Flag(1,4; F)$. It has rank $13$, so the dimension of the desingularization is
$13+15=28$ as needed.
The orbit closure is normal, with rational singularities. The terms in the resulting complex $\FF(6)_\bullet$ are
$$0\rightarrow (7^6,6)\rightarrow (7,6^5,5)\rightarrow (6^2,5^4,4)\rightarrow$$
$$\rightarrow (5^3,4^3,3)\rightarrow (4^4,3^2,2)\rightarrow$$
$$\rightarrow (3^5,2,1)\rightarrow (2^6,0)\rightarrow (0^7).$$

Notice that the complex $\FF(6)_\bullet$ is pure.

\bigskip\bigskip

$\spadesuit$ The orbit $\OOO_{3}$ of $1$-decomposable tensors (codimension $14$).

This orbit closure is the set of tensors $t\in\bigwedge^3\CC^7$ that can be expressed as $t=\ell\wedge\overline t$ where $t\in F$, ${\overline t}\in\bigwedge^2 F$.
The desingularization $Z(3)$ lives on the Grassmannian $Grass(6,F)$. Denoting the tautological bundles as $\Cal R$, $\Cal Q$ ($rank\ {\Cal R}=6$, $rank\ \Cal Q=1$), we have
$\xi =\bigwedge^3{\Cal R}$. The orbit closure has dimension $15+6=21$, so its codimension is $14$. It is normal and has rational singularities.
Calculating the resolution is straightforward, as $\xi$ is irreducible. The defining ideal is generated by the representation $(2^3,1^3,0)$ in degree $3$.

\bigskip\bigskip

$\spadesuit$ The generic degenerate orbit closure $\overline{\Cal O}_5$ of tensors of rank $\le 6$ (codimension $9$).

This orbit closure has a desingularization $Z(5)$ that lives on the Grassmannian $Grass(1,F)$. Denoting the tautological bundles
 $\Cal R$, $\Cal Q$ ($rank\ {\Cal R}=1$, $rank\ \Cal Q=6$), we have
$\xi ={\Cal R}\otimes\bigwedge^2 {\Cal Q}$. The orbit closure has dimension $20+6=26$, so its codimension is $9$. It is normal and has rational singularities.
Calculating the resolution is straightforward, as $\xi$ is irreducible.

\bigskip\bigskip

\head \S 4. The case $(E_7, \alpha_3)$. \endhead

$X=E\otimes\bigwedge^2 F$, $E=\CC^2$,
$F=\CC^6$, ${G}= SL(E)\times SL(F)\times \CC^*$.

The graded Lie algebra of type $E_7$ is
$${\goth g}(E_7)= {\goth g}_{-3}\oplus {\goth g}_{-2}\oplus {\goth g}_{-1}\oplus {\goth g}_0\oplus {\goth g}_1\oplus {\goth g}_2\oplus {\goth g}_3$$
with $G_0=SL(2)\times SL(6)\times\CC^*$, ${\goth g}_0= \CC\oplus {\goth sl}(2)\oplus{\goth sl}(6)$,
${\goth g}_1=\CC^2\otimes\bigwedge^2\CC^6$, ${\goth g}_2=\bigwedge^2\CC^2\otimes\bigwedge^4\CC^6$, ${\goth g}_3=S_{2,1}\CC^2\otimes\bigwedge^6\CC^6$.

Let $\lbrace e_1, e_2\rbrace$ be a basis of $E$, $\lbrace f_1 ,\ldots ,f_6\rbrace$ be a basis of $F$.
We denote the tensor $e_a\otimes f_i\wedge f_j$ by $[a;ij]$.
The invariant scalar product on $\goth h$ restricted to the roots from ${\goth g}_1$ is
$$([a;ij], [b;kl])=\delta-1$$
where $\delta =\#(\lbrace a\rbrace\cap\lbrace b\rbrace )+\#(\lbrace i,j\rbrace\cap\lbrace k,l\rbrace ).$

The ring of invariants is
generated by an invariant (hyperdiscriminant $\Delta$) of degree 12.

This representation has $15$ orbits.

$$\matrix number&{\goth s}&dim&representative\\
0&0&0&0\\
1&A_1&10&[1;12]\\
2&2A_1&15&[1;12]+[1;34] \\
3&2A_1&15&[1;12]+[2;13]\\
4&3A_1&16&[1;12]+[1;34]+[1;56]\\
5&3A_1&19&[1;12]+[1;34]+[2;13]\\
6&A_2&20&[1;12]+[2;34]\\
7&A_2+A_1&23&[1;12]+[2;34]+[1;35].\\
8&A_2+2A_1&25&[1;12]+[2;34]+[1;35]+[2;15]\\
9&A_2+2A_1&24&[1;12]+[2;34]+[1;35]+[1;46]\\
10&2A_2&26&[1;12]+[2;34]+[1;45]+[2;16]\\
11&2A_2+A_1&28&[1;12]+[2;34]+[1;45]+[2;16]+[1;36]\\
12&A_3&25&[1;12]+[2;34]+[1;56]\\
13&A_3+A_1&29&[1;12]+[2;34]+[1;56]+[2;15]\\
14&D_4(a_1)&30&<[1;12],[2;12],[1;34],[2;34],[1;56],[2;56]> \endmatrix$$

Our representation can be treated as a set of skew symmetric $6\times 6$ matrices, with
linear entries in two variables $x,y$. From that point of view the geometry of the orbits is described in the next table.

 $$\matrix number&proj.\ picture&tensor\ picture\\
0&0&0\\
1&C(Seg(\PP^1\times Grass(2,6))&h.w.\ vector\\
2&&$F$-degenerate\\
3&&$F$-degenerate \\
4&&$E$-degenerate\\
5&\tau (\overline\OOO_1)&$F$-degenerate\\
6&\sigma_2 (\overline\OOO_1)&$F$-degenerate&\\
7&&$F$-degenerate\\
8&&$F$-degenerate\\
9&&rank-2-member;Pf=x^3\\
10&&Pf=0\\
11&&Pf=x^3 \\
12&&rank-2-member;Pf=x^2y\\
13&&Pf=x^2y \\
14&\sigma_3 (\overline\OOO_1)&generic  \endmatrix$$

The numerical data is as follows

$$\matrix number&degree&numerator\\
0&1&1\\
1&126&1+20t+60t^2+40t^3+5t^4\\
2&42&1+15t+15t^2+11t^3\\
3&364&1+15t+75t^2+147t^3+105t^4+21t^5\\
4&15&1+14t\\
5&876&1+11t+66t^2+212t^3+316t^4+210t^5+55t^6+5t^7\\
6&633&1+10t+55t^2+146t^3+209t^4+146t^5+55t^6+10t^7+t^8\\
7&588&1+7t+43t^2+113t^3+161t^4+149t^5+82t^6+28t^7+4t^8\\
8&108&1+5t+15t^2+31t^3+35t^4+21t^5\\
9&238&1+6t+36t^2+74t^3+72t^4+42t^5+7t^6\\
10&81&1+4t+10t^2+16t^3+19t^4+16t^5+10t^6+4t^7+t^8\\
11&27&1+2t+3t^2+4t^3+5t^4+6t^5+4t^6+2t^7\\
12&84&1+5t+30t^2+40t^3+8t^4\\
13&12&1+t+t^2+3t^3+3t^4+3t^5\\
14&1&1
\endmatrix$$

The singularities data is

 $$\matrix number&spherical&normal&C-M&R.S.&Gor\\
0&yes&yes&yes&yes&yes\\
1&yes&yes&yes&yes&no\\
2&yes&yes&yes&yes&no\\
3&yes&yes&yes&yes&no\\
4&yes&yes&yes&yes&no\\
5&yes&yes&yes&yes&no\\
6&no&yes&yes&yes&no\\
7&no&no&no&no&no\\
n(7)&no&yes&yes&yes&no\\
8&no&yes&yes&yes&no \\
9&no&no&no&no&no\\
n(9)&no&yes&yes&yes&no\\
10&no&yes&yes&yes&yes\\
11&no&yes&yes&es&no\\
12&no&no&no&no&no\\
n(12)&no&yes&yes&yes&no\\
13&no&no&no&no&yes\\
n(13)&no&yes&yes&yes&no\\
14&yes&yes&yes&yes&yes  \endmatrix$$

\proclaim{Remark}
The degeneration partial order is
$$
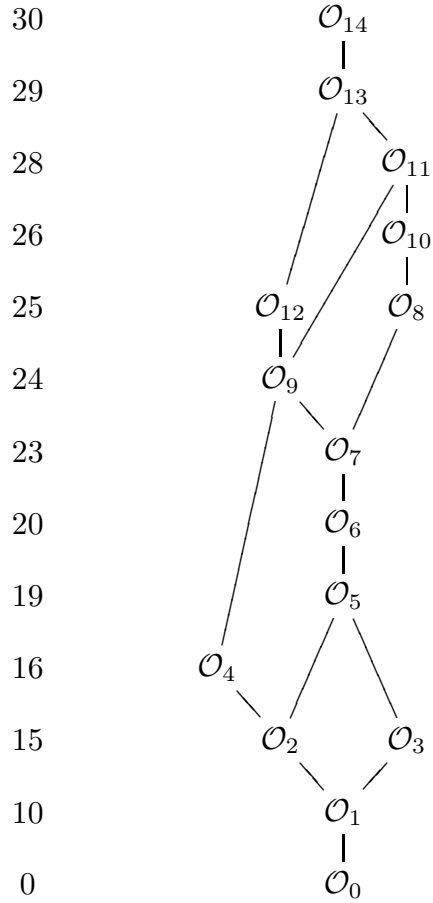
\xy
(20,0)*+{{\Cal O}_{0}}="o0";%
(20,8)*+{{\Cal O}_{1}}="o1";%
(13,16)*+{{\Cal O}_2}="o2";%
(27,16)*+{{\Cal O}_3}="o3";%
(6,24)*+{{\Cal O}_4}="o4";%
(20,32)*+{{\Cal O}_5}="o5";%
(20,40)*+{{\Cal O}_6}="o6";%
(20,48)*+{{\Cal O}_7}="o7";%
(27,64)*+{{\Cal O}_8}="o8";%
(13,56)*+{{\Cal O}_9}="o9";%
(27,72)*+{{\Cal O}_{10}}="o10";%
(27,80)*+{{\Cal O}_{11}}="o11";%
(13,64)*+{{\Cal O}_{12}}="o12";%
(20,88)*+{{\Cal O}_{13}}="o13";%
(20,96)*+{{\Cal O}_{14}}="o14";%
(-15,0)*{0};
(-15,8)*{10};
(-15,16)*{15};
(-15,24)*{16};
(-15,32)*{19};
(-15,40)*{20};
(-15,48)*{23};
(-15,56)*{24};
(-15,64)*{25};
(-15,72)*{26};
(-15,80)*{28};
(-15,88)*{29};
(-15,96)*{30};
{\ar@{-} "o0"; "o1"};%
{\ar@{-} "o1"; "o2"};%
{\ar@{-} "o1"; "o3"};%
{\ar@{-} "o2"; "o4"};%
{\ar@{-} "o2"; "o5"};%
{\ar@{-} "o3"; "o5"};%
{\ar@{-} "o4"; "o9"};%
{\ar@{-} "o5"; "o6"};%
{\ar@{-} "o6"; "o7"};%
{\ar@{-} "o7"; "o8"};%
{\ar@{-} "o7"; "o9"};%
{\ar@{-} "o8"; "o10"};%
{\ar@{-} "o9"; "o12"};%
{\ar@{-} "o9"; "o11"};%
{\ar@{-} "o10"; "o11"};%
{\ar@{-} "o11"; "o13"};%
{\ar@{-} "o12"; "o13"};%
{\ar@{-} "o13"; "o14"};%
\endxy
$$

\endproclaim

$\spadesuit$ The hyperdiscriminant orbit closure $\overline{\OOO_{13}}$.

This is the hypersurface given by the tensors with vanishing hyperdiscriminant.

Its desingularization is, as always (see \cite{KW11a} section 5) is given by the bundle $\eta$ whose complementary bundle is the $1$-jet bundle $\xi(13)$.
The orbit closure is not normal.
The resolution of the normalization is
$$(4,2;2^6)\rightarrow (2,1;1^6)\oplus (0,0;0^6).$$
The defining equation of the orbit itself is the hyperdeterminant which has degree $12$.
Notice that extra partition on the term $\FF(13)_0$ is just ${\goth g}_3$.

\bigskip\bigskip

$\spadesuit$ The codimension $5$ orbit closure $\overline{\OOO_{12}}$.

Take $G/P= \PP (E)\times Grass(4,F)$. Take $\xi(12) = \OOO(-1)\otimes Ker(\bigwedge^2 F\rightarrow \bigwedge^2{\Cal Q})$. The complementary bundle $\eta(12)$ defines desingularization $Z(12)$.
 The terms of the resulting complex $\FF(12)_\bullet$ are
 $$(8,1;3^6)\rightarrow (6,1;3^2,2^4)\rightarrow (5,1;3,2^4,1)\rightarrow$$
 $$\rightarrow (4,1;2^5,0)\oplus (3,1;3,1^5)\rightarrow$$
 $$\rightarrow (2,1;2,1^4,0)\rightarrow (1,1;1^4,0^2)\oplus(0,0;0^6).$$

The orbit closure consists of pencils of skew symmetric matrices such that one of the matrices has rank $\le 2$.
If we treat the tensor from that orbit as a skew symmetric $6\times 6$ matrix of linear forms, the Pfaffian is divisible by $x^2$, but not by $x^3$.
The orbit is not normal, but its normalization has rational singularities.

Using Macaulay 2 (see  \cite{G11}) it is possible to resolve the cokernel $C_{12}$.
It's Betti table (with four rows, fourteen columns)  is
$$\matrix 15&70&.&.&.&.&.\\
.&63&.&.&.&.&.\\
.&84&175&90&.&.&.\\
.&735&4170&11555&21796&32247&38670
\endmatrix$$
$$\matrix
&.&.&.&.&.&.&.\\
&.&.&.&.&.&.&.\\
&.&.&.&.&.&.&.\\
&36295&26452&12858&4556&1110&180&16
\endmatrix$$

Taking mapping cone we conclude that the resolution of $\CC[{\overline {\OOO_{12}}}]$ has only two linear strands,
and the defining ideal of $\overline{\OOO_{12}}$ is generated by 735 equations in degree 5.
\bigskip\bigskip

$\spadesuit$ The codimension $2$ orbit closure $\overline{\OOO_{11}}$.

 The orbit closure ${\overline\OOO}_{11}$ is  the singular locus of ${\overline\OOO}_{13}$. The desingularization $Z(11)$  is given by the bundle $\eta(11)$ whose complementary bundle $\xi(11)$ corresponds to a $P$-module with the following 15 weights:

$$(1,0;0,0,1,0,1,0), (1,0;0,0,0,1,1,0), (1,0;0,0,1,0,0,1),$$
$$ (1,0;0,0,0,1,0,1), (1,0;0,0,0,0,1,1),(0,1;1,0,0,0,1,0), $$
$$ (0,1;1,0,0,0,0,1), (0,1;0,1,0,0,1,0), (0,1;0,1,0,0,0,1), $$
$$(0,1;0,0,1,1,0,0), (0,1;0,0,1,0,1,0), (0,1;0,0,0,1,1,0), $$
$$(0,1;0,0,1,0,0,1), (0,1;0,0,0,1,0,1), (0,1;0,0,0,0,1,1).$$

It lives on the homogeneous space $\PP^1\times Flag(2,4;F)$. So dimension of $Z(11)$ equals $15+1+8+4=28$ as required.

The orbit closure ${\overline\OOO}_{12}$ is normal, with rational singularities.

 The resolution has terms
$$(5,4;3^6)\rightarrow (4,2;2^6)\rightarrow (0,0;0^6).$$
and is determinantal by Hilbert-Burch Theorem. If we treat the tensor from that orbit as a $6\times 6$ skew symmetric matrix of linear forms, the Pfaffian is a  binary cubic. In fact in this resolution all matrix entries are polynomials in the coefficient of the Pfaffian (cubic binary form). The resolution is the same as for the ideal if the set of cubics that are powers of linear form.

\bigskip\bigskip

$\spadesuit$. The codimension $4$ orbit closure $\overline{\Cal O}_{10}$ of tensors with vanishing Pfaffian.

Take $G/P= Grass(4,F)$ and take $\xi(10) =E\otimes\bigwedge^2{\Cal R}$. We get our desingularization $Z(10)$. The complex $\FF(10)_\bullet$ one gets has the terms
$$(6,6;4^6)\rightarrow (6,3;3^6)\rightarrow$$
$$\rightarrow (3,3;2^6)\oplus (5,1;2^6)\rightarrow$$
$$\rightarrow (3,0;1^6)\rightarrow (0,0;0^6).$$

This is a Koszul complex on 4 equations which are the  coefficients of the Pfaffian binary cubic form. The orbit closure is the set of these pencils of skew symmetric $6\times 6$ matrices for which their Pfaffian is identically zero.

There are in fact two orbits with this property such that the pencil does not intersect matrices of rank $\le 2$ (the other is the degenerate orbit $\OOO_9$.). This fact is proved in \cite{MM} .

\bigskip\bigskip

$\spadesuit$ The codimension $6$ orbit closure $\overline{\Cal O}_{9}$.

 Consider the orbit closure ${\overline\OOO}_9$. The desingularization $Z(9)$ is obtained from a bundle $\eta(9)$  with $11$ weights
 $$(1,0;1,1,0,0,0,0), (1,0;1,0,1,0,0,0), (1,0;1,0,0,1,0,0),$$
$$(1,0;1,0,0,0,1,0), (1,0;1,0,0,0,0,1), (1,0;0,1,1,0,0,0),$$
$$(1,0;0,1,0,1,0,0), (1,0;0,1,0,0,1,0), (1,0;0,1,0,0,0,1),$$
$$(1,0;0,0,1,1,0,0), (0,1;1,1,0,0,0,0).$$
This bundle lives on $\PP^1\times Flag(2,4;F)$ so the dimension of the desingularization is $11+1+8+4=24$ as required.

The variety ${\overline\OOO}_9$ is not normal, but the normalization has rational singularities. This follows from Remark 1.1.

Here we can calculate the Euler characteristics of the exterior powers of $\xi(9)$ and of low symmetric powers of $\eta(9)$ only. This proves the orbit closure is not normal.
Based on this one can conjecture the following. The reader will be able to recover the Euler characteristics of the bundles $\bigwedge^i(\xi(9))$ from the data below.

\proclaim{Conjecture}
The terms in the resolution $\FF(9)_\bullet$ are as follows.

$$H^*(\bigwedge^0\xi )=(0,0;0,0,0,0,0,0)[0],$$
$$H^*(\bigwedge^2\xi )=(1,1;1,1,1,1,0,0)[2],$$
$$H^*(\bigwedge^3\xi )=(2,1;1,1,1,1,1,1)[2]\oplus (2,1;2,1,1,1,1,0)[2],$$
$$H^*(\bigwedge^4\xi )=(3,1;3,1,1,1,1,1)[2],$$
$$H^*(\bigwedge^5\xi )=(4,1;2,2,2,2,2,0)[3],$$
$$H^*(\bigwedge^6\xi )=(3,3;3,2,2,2,2,1)[4]\oplus (5,1;3,2,2,2,2,1)[3],$$
$$H^*(\bigwedge^7\xi )=(4,3;4,2,2,2,2,2)[4]\oplus (6,1;3,3,2,2,2,2)[3],$$
$$H^*(\bigwedge^8\xi )=(5,3;3,3,3,3,3,1)[5],$$
$$H^*(\bigwedge^9\xi )=(6,3;4,3,3,3,3,2)[5]\oplus (8,1;3,3,3,3,3,3)[4],$$
$$H^*(\bigwedge^{10}\xi )=(7,3;4,4,3,3,3,3)[5],$$
$$H^*(\bigwedge^{12}\xi )=(9,3;4,4,4,4,4,4)[6].$$

\endproclaim

\bigskip\bigskip

\head \S 5. The case $(E_7, \alpha_4)$. \endhead

  $X= E\otimes F\otimes H$, $E=\CC^2 ,F=\CC^3
,H=\CC^4$, ${G_0}=SL(E)\times SL(F)\times SL(H)\times \CC^*$.

The graded Lie algebra of type $E_7$ is
$${\goth g}(E_7)= {\goth g}_{-4}\oplus {\goth g}_{-3}\oplus {\goth g}_{-2}\oplus {\goth g}_{-1}\oplus {\goth g}_0\oplus {\goth g}_1\oplus {\goth g}_2\oplus {\goth g}_3\oplus{\goth g}_4$$
with ${\goth g}_0= \CC\oplus {\goth sl}(2)\oplus{\goth sl}(3)\oplus{\goth sl}(4)$,
${\goth g}_1=\CC^2\otimes\CC^3\otimes\CC^4$, ${\goth g}_2=\bigwedge^2\CC^2\otimes\bigwedge^2\CC^3\otimes\bigwedge^2\CC^4$, ${\goth g}_3=S_{2,1}\CC^2\otimes\bigwedge^3\CC^3\otimes\bigwedge^3\CC^4$, ${\goth g}_4 =S_{2,2}\CC^2\otimes S_{2,1,1}\CC^3\otimes\bigwedge^4\CC^4$.

Let $\lbrace e_1, e_2\rbrace$ be a basis of $E$, and $\lbrace f_1, f_2, f_3\rbrace$,  $\lbrace h_1, h_2, h_3, h_4\rbrace$ bases of $F$, $H$ respectively.
We label $e_a\otimes f_i\otimes h_u$ by $[a;i;u]$.

The invariant scalar product on $\goth h$ restricted to the roots from ${\goth g}_1$ is
$$([a;i;u], [b;j;v])=\delta-1$$
where $\delta =\#(\lbrace a\rbrace\cap\lbrace b\rbrace )+\#(\lbrace i\rbrace\cap\lbrace j\rbrace )+\#(\lbrace u\rbrace\cap\lbrace v\rbrace ).$

There are six $H$-nondegenerate orbits. They can be described by observing that the
castling transform establishes a bijection between $H$-nondegenerate orbits and
$H'$-nondegenerate orbits for the $2\times 3\times 2$ matrices corresponding to representation $E\otimes F\otimes H'$.
The six orbits in the representation $E\otimes F\otimes H'$ are: generic, hyperdeterminant hypersurface and four $F$-degenerate orbits, coming from $2\times 2\times 2$ matrices: generic, hyperdeterminant and two determinantal varieties.
Combining this knowledge with the case $(E_6,4)$ we get  $24$ orbits in our representation.

$$\matrix number&\goth s &dim&representative\\
0&0&0&0\\
1&A_1&7&[1;1;1]\\
2&2A_1&9&[1;1;1]+[2;2;1]\\
3&2A_1&10&[1;1;1]+[2;1;2]\\
4&2A_1&11&[1;1;1]+[1;2;2]\\
5&3A_1&13&[1;1;1]+[1;2;2]+[2;1;2]\\
6&3A_1&13&[1;1;1]+[1;2;2]+[1;3;3]\\
7&A_2&14&[1;1;1]+[2;2;2]\\
8&A_2+A_1&15&[1;1;1]+[2;2;2]+[1;3;2]\\
9&A_2+A_1&16&[1;1;1]+[2;2;2]+[1;2;3] \\
10&A_2+2A_1&16&[1;1;1]+[2;2;2]+[1;3;2]+[2;3;1] \\
11&A_2+2A_1&17&[1;1;1]+[2;2;2]+[1;2;3]+[2;1;3] \\
12&A_2+2A_1&17&[1;1;1]+[2;2;2]+[1;2;3]+[1;3;2]\\
13&2A_2&17&[1;1;1]+[2;2;2]+[1;2;3]+[2;3;1]\\
14&A_3&18&[1;1;1]+[2;2;2]+[1;3;3]\\
15&2A_2+A_1&19&[1;1;1]+[2;2;2]+[1;2;3]+[2;3;1]+[1;3;2]\\
16&A_3+A_1&20&[1;1;1]+[2;2;2]+[1;3;3]+[2;1;3]\\
17&D_4(a_1)&21&[1;1;1]+[2;1;1]+[1;2;2]-[2;2;2]+[1;3;3] \\
18&2A_2&18&[1;1;1]+[2;2;2]+[2;1;3]+[1;2;4]\\
19&A_3+A_1&19&[1;1;1]+[2;2;2]+[1;3;3]+[1;2;4]\\
20&A_3+A_2&21&[1;1;1]+[2;2;2]+[1;3;3]+[2;1;3]+[1;2;4]\\
21&A_4&22&[1;1;1]+[2;2;2]+[1;3;3]+[2;1;4]\\
22&A_4+A_1&23&[1;1;1]+[2;2;2]+[1;3;3]+[2;1;4]+[1;2;4]\\
23&A_4+A_2&24&[1;1;1]+[2;2;2]+[1;3;3]+[2;1;4]+[1;2;4]\\
&&&+[2;3;1]  \endmatrix$$

Our representation has two interpretations. One can view it as a space of $3\times 4$ matrices with the entries
that are linear forms in two variables $x, y$, and as a space of quiver representations of Kronecker quiver of dimenson vector $(3,4)$.
In describing the orbits we refer to ''a matrix picture" and ''a quiver picture" to refer to these interpretations.

$$\matrix number&proj.\ picture&matrix\ pic.&quiver\ pic.\\
0&0&&\\
1&C(Seg(\PP^1\times \PP^2\times\PP^3))&\overline {\OOO}_1 ,(E_6 ,4)&\\
2&&\overline {\OOO}_2 ,(E_6 ,4)&\\
3&&\overline {\OOO}_3 ,(E_6 ,4)&\\
4&&\overline {\OOO}_4 ,(E_6 ,4)&\\
5&\tau(\overline{\OOO_1})&\overline {\OOO}_5 ,(E_6 ,4)&\\
6&&\overline {\OOO}_6 ,(E_6 ,4)&\\
7&\sigma_2(\overline{\OOO_1})&\overline {\OOO}_7 ,(E_6 ,4)&\\
8&&\overline {\OOO}_8 ,(E_6 ,4)&\\
9&&\overline {\OOO}_9 ,(E_6 ,4)& \\
10&&\overline {\OOO}_{10} ,(E_6 ,4) \\
11&&\overline {\OOO}_{11} ,(E_6 ,4) \\
12&&\overline {\OOO}_{12} ,(E_6 ,4)&\\
13&&\overline {\OOO}_{13} ,(E_6 ,4)&\\
14&&\overline {\OOO}_{14} ,(E_6 ,4)&\\
15&&\overline {\OOO}_{15} ,(E_6 ,4)&\\
16&&\overline {\OOO}_{16} ,(E_6 ,4)&\\
17&\sigma_3(\overline{\OOO_1})&\overline {\OOO}_{17} ,(E_6 ,4)&\\
18&&$F$-degenerate&(1,0)\oplus (2,4)\\
19&&rank-1-member&\\
20&&rank-2-member&(2,2)\oplus (1,2)\\
21&&&2*(1,1)\oplus (1,2)\\
22&&hyperdisc.&(1,1)\oplus (2,3)\\
23&\sigma_4(\overline{\OOO_1})&generic&
  \endmatrix$$

  $$\matrix number&tensor\ picture\\
0&0\\
1&h.w.\ vector\\
2&$H$-rank\le 1\\
3&$H$-rank\le 3, $F$-rank\le 1\\
4&$H$-rank\le 3, $E$-rank\le 1, $F-H$-rank\le 2\\
5&$H$-rank\le 3, hyperdet\ in\ \OOO_6 \\
6&$E$-rank\le 1&\\
7&\CC^2\otimes\CC^2\otimes\CC^2\\
8&$H$-rank\le 3, hyperdet\ in\ \OOO_{10}\\
9&$H$-rank\le 3, hyperdet\ in \OOO_{11} \\
10&$H$-rank\le 2 \\
11&$H$-rank\le 3, $F$-rank\le 2 \\
12&\overline {\OOO}_{12} ,(E_6 ,3)\\
13&\overline {\OOO}_{13} ,(E_6 ,3)\\
14&\overline {\OOO}_{14} ,(E_6 ,3)\\
15&\overline {\OOO}_{15} ,(E_6 ,3)\\
16&$H$-rank\le 3, hyperdisc.\ \CC^2\otimes\CC^3\otimes\CC^3\\
17&$H$-rank\le 3, generic\ \CC^2\otimes\CC^3\otimes\CC^3\\
18&$F$-rank\le 2\\
19&rank-1\ member\\
20&sing.locus\ \PP^1\\
21&sing(hyperdisc.)\\
22&hyperdisc\\
23&generic  \endmatrix$$

The numerical data are as follows.

$$\matrix number&degree&numerator\\
0&1&1\\
1&60&1+17t+33t^2+9t^3\\
2&56&1+15t+30t^2+10t^3\\
3&36&1+14t+21t^2\\
4&60&1+13t+25t^2+21t^3\\
5&408&1+11t+66t^2+166t^3+131t^4+33t^5\\
6&12&1+11t\\
7&276&1+10t+55t^2+100t^3+85t^4+22t^5+3t^6\\
8&420&1+9t+63t^2+135t^3+150t^4+54t^5+8t^6\\
9&240&1+8t+54t^2+96t^3+81t^4\\
10&105&1+8t+36t^2+40t^3+20t^4\\
11&96&1+7t+28t^2+28t^3+22t^4+10t^5\\
12&312&1+7t+46t^2+98t^3+113t^4+47t^5\\
13&216&1+7t+46t^2+82t^3+70t^4+10t^5\\
14&192&1+6t+39t^2+60t^3+51t^4+30t^5+5t^6\\
15&288&1+5t+15t^2+35t^3+55t^4+75t^5+65t^6+37t^7\\
16&180&1+4t+10t^2+28t^3+40t^4+52t^5+30t^6+12t^7+3t^8\\
17&20&1+3t+6t^2+10t^3\\
18&28&1+6t+21t^2\\
19&60&1+5t+33t^2+21t^3\\
20&88&1+3t+6t^2+18t^3+27t^4+33t^5\\
21&36&1+2t+3t^2+12t^3+9t^4+6t^5+3t^6\\
22&12&1+t+t^2+9t^3\\
23&1&1 \endmatrix$$

The singularities data are as follows.

  $$\matrix number&spherical&normal&C-M&R.S.&Gor\\
0&yes&yes&yes&yes&yes\\
1&yes&yes&yes&yes&no\\
2&yes&yes&yes&yes&no\\
3&yes&yes&yes&yes&no\\
4&yes&yes&yes&yes&no\\
5&yes&yes&yes&yes&no\\
6&yes&yes&yes&yes&no\\
7&no&yes&yes&yes&no\\
8&no&no&no&no&no\\
n(8)&no&yes&yes&yes&no\\
9&no&no&no&no&no \\
n(9)&no&yes&yes&yes&no\\
10&no&yes&yes&yes&no \\
11&no&yes&yes&yes&no \\
12&no&no&no&no&no\\
n(12)&no&yes&yes&yes&no\\
13&no&no&no&no&no\\
n(13)&no&yes&yes&yes&no\\
14&no&no&no&no&no\\
n(14)&no&yes&yes&yes&no\\
15&no&yes&yes&yes&no\\
16&no&no&no&no&no\\
n(16)&no&yes&yes&yes&no\\
17&no&yes&yes&yes&no\\
 18&no&yes&yes&yes&no\\
19&no&no&no&no&no\\
n(19)&no&yes&yes&yes&no\\
20&no&no&no&no&no\\
n(20)&no&yes&yes&yes&no\\
21&no&no&yes&no&no\\
n(21)&no&yes&yes&yes&no\\
22&no&no&yes&no&yes\\
n(22)&no&yes&yes&yes&no\\
23&no&yes&yes&yes&yes \endmatrix$$

\proclaim{Remark} The degeneration order is
$$
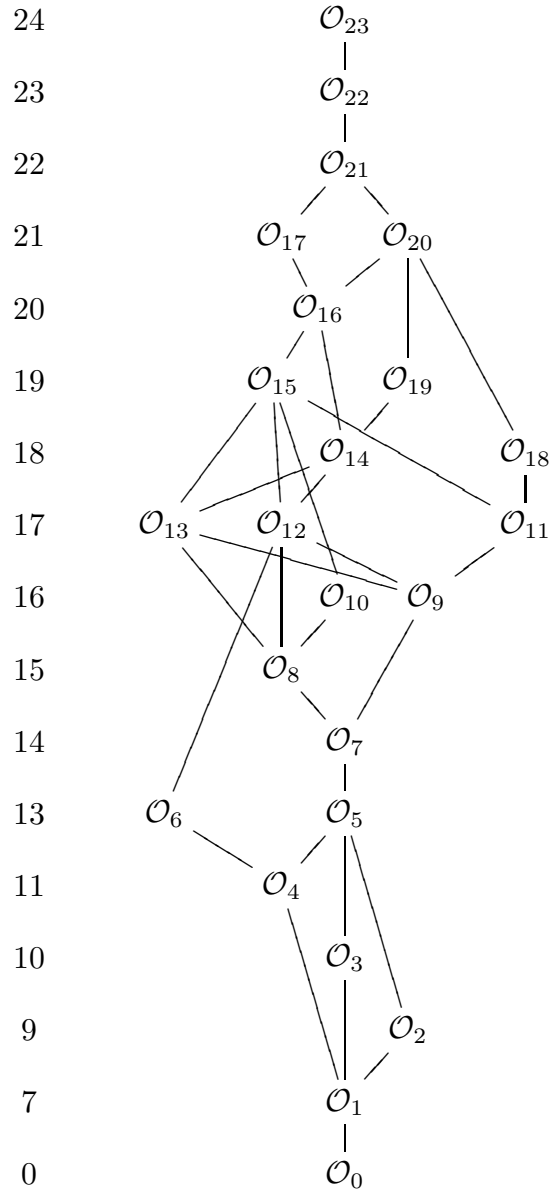
\xy
(20,0)*+{{\Cal O}_{0}}="o0";%
(20,8)*+{{\Cal O}_{1}}="o1";%
(27,16)*+{{\Cal O}_2}="o2";%
(20,24)*+{{\Cal O}_3}="o3";%
(13,32)*+{{\Cal O}_4}="o4";%
(20,40)*+{{\Cal O}_5}="o5";%
(0,40)*+{{\Cal O}_6}="o6";%
(20,48)*+{{\Cal O}_7}="o7";%
(13,56)*+{{\Cal O}_8}="o8";%
(29,64)*+{{\Cal O}_9}="o9";%
(20,64)*+{{\Cal O}_{10}}="o10";%
(40,72)*+{{\Cal O}_{11}}="o11";%
(13,72)*+{{\Cal O}_{12}}="o12";%
(0,72)*+{{\Cal O}_{13}}="o13";%
(20,80)*+{{\Cal O}_{14}}="o14";%
(12,88)*+{{\Cal O}_{15}}="o15";%
(17,96)*+{{\Cal O}_{16}}="o16";%
(13,104)*+{{\Cal O}_{17}}="o17";%
(40,80)*+{{\Cal O}_{18}}="o18";%
(27,88)*+{{\Cal O}_{19}}="o19";%
(27,104)*+{{\Cal O}_{20}}="o20";%
(20,112)*+{{\Cal O}_{21}}="o21";%
(20,120)*+{{\Cal O}_{22}}="o22";%
(20,128)*+{{\Cal O}_{23}}="o23";%
(-15,0)*{0};
(-15,8)*{7};
(-15,16)*{9};
(-15,24)*{10};
(-15,32)*{11};
(-15,40)*{13};
(-15,48)*{14};
(-15,56)*{15};
(-15,64)*{16};
(-15,72)*{17};
(-15,80)*{18};
(-15,88)*{19};
(-15,96)*{20};
(-15,104)*{21};
(-15,112)*{22};
(-15,120)*{23};
(-15,128)*{24};
{\ar@{-} "o0"; "o1"};%
{\ar@{-} "o1"; "o2"};%
{\ar@{-} "o1"; "o3"};%
{\ar@{-} "o1"; "o4"};%
{\ar@{-} "o2"; "o5"};%
{\ar@{-} "o3"; "o5"};%
{\ar@{-} "o4"; "o5"};%
{\ar@{-} "o4"; "o6"};%
{\ar@{-} "o5"; "o7"};%
{\ar@{-} "o6"; "o12"};%
{\ar@{-} "o7"; "o8"};%
{\ar@{-} "o7"; "o9"};%
{\ar@{-} "o8"; "o10"};%
{\ar@{-} "o8"; "o12"};%
{\ar@{-} "o8"; "o13"};%
{\ar@{-} "o9"; "o11"};%
{\ar@{-} "o9"; "o12"};%
{\ar@{-} "o9"; "o13"};%
{\ar@{-} "o10"; "o15"};%
{\ar@{-} "o11"; "o15"};%
{\ar@{-} "o11"; "o18"};%
{\ar@{-} "o12"; "o14"};%
{\ar@{-} "o12"; "o15"};%
{\ar@{-} "o13"; "o14"};%
{\ar@{-} "o13"; "o15"};%
{\ar@{-} "o14"; "o16"};%
{\ar@{-} "o14"; "o19"};%
{\ar@{-} "o15"; "o16"};%
{\ar@{-} "o16"; "o17"};%
{\ar@{-} "o16"; "o20"};%
{\ar@{-} "o17"; "o21"};%
{\ar@{-} "o18"; "o20"};%
{\ar@{-} "o19"; "o20"};%
{\ar@{-} "o20"; "o21"};%
{\ar@{-} "o21"; "o22"};%
{\ar@{-} "o22"; "o23"};%
\endxy
$$
\endproclaim

\bigskip\bigskip

Next we will describe in detail the non-degenerate orbit closures in $\CC^2\otimes\CC^3\otimes\CC^4$.
These are the orbits $\overline{\OOO_{18}}$, $\overline{\OOO_{19}}$. $\overline{\OOO_{20}}$, $\overline{\OOO_{21}}$ and $\overline{\OOO_{22}}$.
The orbit $\OOO_{23}$ is generic so there is not much to say.

We use the usual notation. $A=Sym(E^*\otimes F^*\otimes H^*)$ and $(a,b;c,d,e,f,g,h,i)$ abbreviates for $S_{a,b}E^*\otimes S_{c,d,e}F^*\otimes S_{f,g,h,i}H^*$.

$\spadesuit$ The hyperdiscriminant orbit closure $\overline{\OOO_{22}}$.

This is the hypersurface given by the tensors with vanishing hyperdiscriminant.
Its desingularization is, as always (see \cite{KW11a} section 5), given by the bundle $\eta(22)$ which is complementary to the $1$-jet bundle $\xi(22)$.
The orbit closure is not normal
but its normalization has rational singularities.

The minimal resolution $\FF(22)_\bullet$ of the normalization has terms
$$0\rightarrow (3,1;2,1,1;1,1,1,1))\rightarrow
(2,1;1,1,1;1,1,1,0)\oplus (0,0;0,0,0;0,0,0,0).$$
The extra representation in $\FF(22)_\bullet$ is just ${\goth g}_3$.

\bigskip\bigskip

$\spadesuit$ The codimension $2$ orbit closure $\overline{\OOO_{21}}$.

The orbit closure ${\overline\OOO}_{21}=Sing ({\overline\OOO}_{22})$ which is, in the classification of
\cite{WZ96} denoted $X_{node}(\emptyset )$. It is an orbit closure of codimension 2. The minimal
resolution of the coordinate ring is
$$0\rightarrow (5,4;3,3,3;3,2,2,2)\rightarrow
 (5,3;3,3,2;2,2,2,2)\rightarrow (0,0;0,0,0;0,0,0,0).$$
It is a determinantal ring.

The orbit closure has a desingularization. It comes from a bundle $\xi = E^*\otimes {\Cal R}\otimes {\Cal R}'$
on $Grass(2, F^*)\times Grass(2, H^*)$.
The complex resolving the coordinate ring of the normalization of $\overline{\OOO_{20}}$ is

$$0\rightarrow (4,4;3,3,2;2,2,2,2)\rightarrow (2,2;2,1,1;1,1,1,1)\oplus (3,1;2,1,1;1,1,1,1)\rightarrow $$
$$\rightarrow (2,1;1,1,1;1,1,1,0)\oplus (0^2;0^3;0^4).$$

\bigskip\bigskip

$\spadesuit$ The codimension $3$ orbit closure $\overline{\OOO_{20}}$.

The orbit closure ${\overline\OOO}_{20}$ of tensors $\phi$ whose singular locus is of the form $pt\times pt\times K^2$.

The desingularization is given by the vector bundle $\xi(20)$ which is a sum of $1$-jet bundles at points $(0:1)(0:0:1)(0:0:0:1)$ and $(0:1)(0:0:1)(0:0:1:0)$,
thus having the weights; $(1,0;0,0,1;0, 0, 1, 0), (1,0;0,0,1;0,0,0,1)$ and eight weights
$$(0,1;1,0,0;0,0,1,0), (0,1;1,0,0;0,0,0,1), (0,1;0,1,0;0,0,1,0), $$
$$(0,1;0,1,0;0,0,0,1), (0,1;0,0,1;1,0,0,0), (0,1;0,0,1;0,1,0,0), $$
$$(0,1;0,0,1;0,0,1,0), (0,1;0,0,1;0,0,0,1).$$
The terms of the complex are
$$H^* (\bigwedge^0\xi )=(0,0;0,0,0;0,0,0,0)[0],$$
$$H^* (\bigwedge^3\xi )=(2,1;1,1,1;1,1,1,0)[3],$$
$$H^* (\bigwedge^4\xi )=(2,2;2,1,1;1,1,1,1)[3]\oplus (3,1;2,1,1;1,1,1,1)[3],$$
$$H^* (\bigwedge^6\xi )=(3,3;2,2,2;2,2,1,1)[5]\oplus (4,2;2,2,2;2,2,2,0)[5],$$
$$H^* (\bigwedge^7\xi )=(4,3;3,2,2;2,2,2,1)[5]\oplus (5,2;3,2,2;2,2,2,1)[5],$$
$$H^* (\bigwedge^8\xi )=(5,3;4,2,2;2,2,2,2)[5]\oplus (6,2;3,3,2;2,2,2,2)[5].$$

The orbit is not normal, but its normalization has rational singularities.

\bigskip\bigskip

$\spadesuit$  The codimension $5$ orbit closure $\overline{\Cal O}_{19}$ of pencils of $3\times 4$ matrices with a member of rank $1$.

The desingularization of ${\overline\OOO}_{19}$ is obtained from $G/P= Grass(1,E)\times Grass(1 ,F)$. The bundle $\xi(19)$ is just $\OOO(-1)\otimes{\Cal Q}^*\otimes H^*$. The orbit closure is not normal, but the normalization has rational singularities. The terms of the resolution are

$$H^* (\bigwedge^0\xi )=(0,0;0,0,0;0,0,0,0)[0],$$
$$H^* (\bigwedge^2\xi )=(1,1;1,1,0;1,1,0,0)[2],$$
$$H^* (\bigwedge^3\xi )=(2,1;2,1,0;1,1,1,0)[2]\oplus (2,1;1,1,1;2,1,0,0)[2],$$
$$H^* (\bigwedge^4\xi )=(3,1;2,1,1;2,1,1,0)[2]\oplus (3,1;3,1,0;1,1,1,1)[2],$$
$$H^* (\bigwedge^5\xi )=(4,1;3,1,1;2,1,1,1)[2],$$
$$H^* (\bigwedge^6\xi )=(5,1;2,2,2;2,2,2,0)[3],$$
$$H^* (\bigwedge^7\xi )=(6,1;3,2,2;2,2,2,1)[3],$$
$$H^* (\bigwedge^8\xi )=(7,1;3,3,2;2,2,2,2)[3],$$

\bigskip\bigskip

$\spadesuit$. The codimension $6$ orbit closure $\overline{\Cal O}_{18}$.

The orbit closure ${\overline\OOO}_{18}$. It is a determinantal variety
of tensors $\phi$ for which the tensor ${\tilde\phi}_{3,1}:E\otimes H\rightarrow F^*$ has rank $\le 2$.
This is a determinantal variety of codimension 6. It is normal, with rational singularities.

\bigskip\bigskip

\head \S 6. The case $(E_7, \alpha_5)$. \endhead

 $X=E\otimes\bigwedge^2 F$, $E=\CC^3$,
$F=\CC^5$, ${ G}= SL(E)\times SL(F)\times \CC^*$.

The graded Lie algebra of type $E_7$ is
$${\goth g}(E_7)= {\goth g}_{-3}\oplus {\goth g}_{-2}\oplus {\goth g}_{-1}\oplus {\goth g}_0\oplus {\goth g}_1\oplus {\goth g}_2\oplus {\goth g}_3$$
with $G_0= SL(3)\times SL(5)\times\CC^*$, ${\goth g}_0= {\goth sl}(3)\oplus{\goth sl}(5)\oplus\CC$,
${\goth g}_1=\CC^3\otimes\bigwedge^2\CC^5$, ${\goth g}_2=\bigwedge^2\CC^3\otimes\bigwedge^4\CC^5$, ${\goth g}_3=\bigwedge^3\CC^3\otimes S_{2,1^4}\CC^5$.

Let $\lbrace e_1, e_2, e_3\rbrace$ be a basis of $E$, $\lbrace f_1 ,\ldots , f_5\rbrace$ be a basis of $F$.
We denote the tensor $e_a\otimes f_i\wedge f_j$ by $[a;ij]$.
The invariant scalar product on $\goth g$ restricted to ${\goth g}_1$ is
$$([a;ij], [b;kl])=\delta-1$$
where $\delta =\#(\lbrace a\rbrace\cap\lbrace b\rbrace )+\#(\lbrace i,j\rbrace\cap\lbrace k,l\rbrace ).$

The non-degenerate orbits for this action were classified by Eisenbud and Koh \cite{EK}. An interesting
feature is that the orbit classification depends on the characteristic of the base
field. Here we just work over $\CC$, so we need the characteristic zero part of the description.
The description of the degenerate orbits follows from the earlier cases.
So we have 8 $E$-degenerate orbits (including zero) coming from the case $E_6, k=3,5$.
We label them ${\overline\OOO}_i (E_6,\alpha_3 )$.
We also have 5 $F$-degenerate orbits which are not $E$-degenerate, coming from the case $(D_6 ,\alpha_3 )$.

 They come from the orbits $O(3,3)$, $O(3,2)$, $O(3,1)$ and $O(3,0)^{\pm}$, and we label them by $\overline{\OOO}(i,j)(D_6,\alpha_3 )$.

 $$\matrix number&{\goth s}&dim&representative\\
0&&0&0\\
1&A_1&9&[1;12]\\
2&2A_1&12&[1;12]+[1;34]\\
3&2A_1&14&[1;12]+[2;13]\\
4&3A_1&15&[1;12]+[2;13]+[3;23]\\
5&3A_1&16&[1;12]+[2;13]+[3;14]\\
6&3A_1&17&[1;12]+[1;34]+[2;13]\\
7&A_2&18&[1;12]+[2;34]\\
8&4A_1&19&[1;12]+[1;34]+[2;13]+[3;14] \\
9&A_2+A_1&20&[1;12]+[2;34]+[1;35]\\
10&A_2+A_1&21&[1;12]+[2;34]+[3;13]\\
11&A_2+2A_1&22&[1;12]+[2;34]+[1;35]+[2;15]\\
12&A_2+2A_1&22&[1;12]+[2;34]+[3;13]+[3;24]\\
13&A_2+2A_1&22&[1;12]+[2;34]+[3;13]+[1;35]\\
14&2A_2&23&[1;12]+[2;34]+[1;45]+[3;13]\\
15&A_3&23&[1;12]+[2;34]+[3;15]\\
16&A_2+3A_1&23&[1;12]+[2;34]+[3;13]+[1;35]+[2;15]\\
17&2A_2+A_1&25&[1;12]+[2;34]+[1;35]+[3;24]+[2;25]\\
18&A_3+A_1&24&[1;12]+[2;34]+[3;15]+[2;25]\\
19&A_3+A_1&26&[1;12]+[2;34]+[3;15]+[1;35]\\
20&A_3+2A_1&27&[1;12]+[2;34]+[3;15]+[1;35]+[2;25]\\
21&D_4(a_1)&27&<[1;12],[1;13],[2;24],[2,34],[3;25],[3;35]>\\
22&D_4 (a_1)+A_1&28&<[1;12],[1;13],[2;24],[2,34],[3;25],[3;35]>+\\
&&&+[1;45]\\
23&A_3+A_2&29&[1;12]+[2;34]+[3;15]+[1;35]+[3;24]\\
24&A_3+A_2+A_1&30&[1;12]+[2;34]+[3;15]+[1;35]+[3;24]+\\
&&&+[2;25]
\endmatrix$$

\proclaim {Remark}
The non-degenerate orbits were classified also by Eisenbud and Koh \cite{EK94} Their classification was based on the possibilities for the scheme $Y$ of $4\times 4$ Pfaffians of a $5\times 5$ skew-symmetric matrix of linear forms. They divided to the following types:
\item{I.} $Y=\emptyset$,
\item{II.a)} $Y_{red}$ is a point ($4$ orbits),
\item{II.b)} $Y_{red}$ is a pair of points ($2$ orbits),
\item{III.} $Y_{red}$ is a line ($3$ orbits),
\item{IV.} $Y_{red}$ spans a plane ($2$ orbits).

In the following table we indicate the Eisendbud-Koh types of non-degenerate orbits.
\endproclaim

The following table gives also geometric descriptions.

Notice that the actual representative of $\overline{\OOO_{21}}$ is $[1;12]+[2;34]+[3;25]+[3;35]$
which is a sum of three decomposable tensor and thus the generic element in $\sigma_3 ({\overline{\OOO_1}})$.

$$\matrix number&proj.\ picture&description&matrix\ picture\\
0&0&{\overline\OOO}_0(E_6,\alpha_3 )&\\
1&h.weight\ vector&{\overline\OOO}_1(E_6,\alpha_3 )&&\\
&C(Seg(\PP^2\times Grass(2,5)))&&&\\
2&&{\overline\OOO}_2(E_6,\alpha_3 )&\\
3&&{\overline\OOO}_3(E_6,\alpha_3 )&\\
4&&{\overline\OOO}(3,0)^+(D_6,\alpha_3 )&\\
5&&{\overline\OOO}(3,0)^-(D_6,\alpha_3 )&\\
6&\tau ({\overline{\OOO_1}})&{\overline\OOO}_4(E_6,\alpha_3 )&\\
7&\sigma_2 ({\overline{\OOO_1}})&{\overline\OOO}_5(E_6,\alpha_3 )&\\
8&&{\overline\OOO}(3,1)(D_6,\alpha_3 )&\\
9&&{\overline\OOO}_6(E_6,\alpha_3 )&\\
10&&{\overline\OOO}(3,2)(D_6,\alpha_3 )&\\
11&\sigma_2 (\overline\OOO_2)&{\overline\OOO}_7(E_6,\alpha_3 )&\\
12&&{\overline\OOO}(3,3)(D_6,\alpha_3 )&\\
13&&III.3.&\\
14&&III.2.&\\
15&&IV.1.&\\
16&&IIa)1.&\\
17&&IIa)2.&\\
18&&III.1.&\\
19&&IIb)1.&\\
20&&IIa)3.&\\
21&\sigma_3 ({\overline{\OOO_1}})&IV.2.&\\
22&&IIb)2.&\\
23&&IIa)4&\\
24&&I&
\endmatrix$$

\bigskip

The numerical data is given in the following table

$$\matrix number&degree&numerator\\
0&1&1\\
1&140&1+21t+66t^2+46t^3+6t^4\\
2&55&1+18t+36t^2\\
3&780&1+16t+106t^2+266t^3+266t^4+110t^5+15t^6\\
4&596&1+15t+90t^2+210t^3+195t^4+75t^5+10t^6\\
5&300&1+14t+75t^2+140t^3+70t^4\\
6&1440&1+13t+91t^2+295t^3+490t^4+400t^5+135t^6+15t^7\\
7&810&1+12t+78t^2+204t^3+285t^4+180t^5+50t^6\\
8&2100&1+11t+66t^2+246t^3+531t^4+633t^5+418t^6+\\
&&+158t^7+33t^8+3t^9\\
9&525&1+10t+70t^2+160t^3+220t^4+64t^5\\
10&855&1+9t+45t^2+125t^3+195t^4+195t^5+160t^6+\\
&&+90t^7+30t^8+5t^9\\
11&45&1+8t+36t^2\\
12&195&1+8t+36t^2+80t^3+70t^4\\
13&2250&1+8t+51t^2+195t^3+435t^4+621t^5+543t^6+\\
&&+291t^7+90t^8+15t^9\\
14&1065&1+7t+43t^2+144t^3+270t^4+315t^5+205t^6+\\
&&+70t^7+10t^8\\
15&1020&1+7t+43t^2+149t^3+275t^4+320t^5+180t^6+45t^7\\
16&825&1+7t+28t^2+84t^3+180t^4+255t^5+190t^6+\\
&&+70t^7+10t^8\\
17&729&1+5t+15t^2+35t^3+70t^4+120t^5+165t^6+\\
&&+165t^7+105t^8+40t^9+8t^{10}\\
18&505&1+6t+36t^2+106t^3+156t^4+150t^5+50t^6\\
19&480&1+4t+10t^2+25t^3+55t^4+94t^5+121t^6+\\
&&+100t^7+55t^8+15t^9\\
20&160&1+3t+6t^2+15t^3+30t^4+45t^5+45t^6+15t^7\\
21&60&1+3t+6t^2+10t^3+15t^4+15t^5+10t^6\\
22&60&1+2t+3t^2+9t^3+15t^4+15t^5+15t^6\\
23&15&1+t+t^2+6t^3+6t^4\\
24&1&1
\endmatrix$$

The singularities data are given in the following table.

  $$\matrix number&spherical&normal&C-M&R.S.&Gor\\
0&yes&yes&yes&yes&yes\\
1&yes&yes&yes&yes&no\\
2&yes&yes&yes&yes&no\\
3&yes&yes&yes&yes&no\\
4&yes&yes&yes&yes&no\\
5&yes&yes&yes&yes&no\\
6&yes&yes&yes&yes&no\\
7&no&yes&yes&yes&no\\
8&yes&yes&yes&yes&no\\
9&no&no&no&no&no \\
n(9)&no&yes&yes&yes&no\\
10&no&yes&yes&yes&no \\
11&no&yes&yes&yes&no \\
12&no&yes&yes&yes&no\\
13&no&no&no&no&no\\
n(13)&no&yes&yes&yes&no\\
14&no&no&no&no&no\\
n(14)&no&yes&yes&yes&no\\
15&no&no&no&no&no\\
n(15)&no&yes&yes&yes&no\\
16&no&yes&yes&yes&no\\
17&no&yes&yes&yes&no\\
 18&no&no&no&no&no\\
 n(18)&no&yes&yes&yes&no\\
19&no&no&no&no&no\\
n(19)&no&yes&yes&yes&no\\
20&no&no&no&no&no\\
n(20)&no&yes&yes&yes&no\\
21&no&yes&yes&yes&no\\
22&no&no&yes&no&no\\
n(22)&no&yes&yes&yes&no\\
23&no&no&yes&yes&yes\\
n(23)&no&yes&yes&yes&no\\
24&no&yes&yes&yes&yes
 \endmatrix$$

 \proclaim{Remark}
The degeneration partial order is

$$
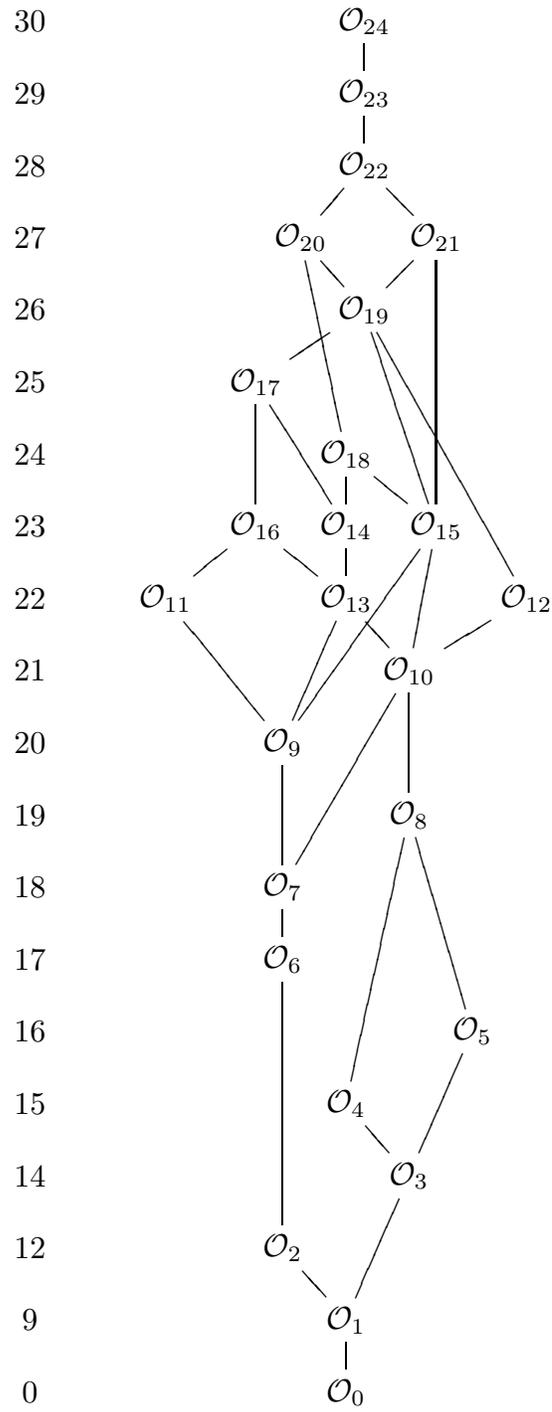
\xy
(20,0)*+{{\Cal O}_{0}}="o0";%
(20,8)*+{{\Cal O}_{1}}="o1";%
(13,16)*+{{\Cal O}_2}="o2";%
(27,24)*+{{\Cal O}_3}="o3";%
(20,32)*+{{\Cal O}_4}="o4";%
(34,40)*+{{\Cal O}_5}="o5";%
(13,48)*+{{\Cal O}_6}="o6";%
(13,56)*+{{\Cal O}_7}="o7";%
(27,64)*+{{\Cal O}_8}="o8";%
(13,72)*+{{\Cal O}_9}="o9";%
(27,80)*+{{\Cal O}_{10}}="o10";%
(0,88)*+{{\Cal O}_{11}}="o11";%
(40,88)*+{{\Cal O}_{12}}="o12";%
(20,88)*+{{\Cal O}_{13}}="o13";%
(20,96)*+{{\Cal O}_{14}}="o14";%
(30,96)*+{{\Cal O}_{15}}="o15";%
(10,96)*+{{\Cal O}_{16}}="o16";%
(10,112)*+{{\Cal O}_{17}}="o17";%
(20,104)*+{{\Cal O}_{18}}="o18";%
(22,120)*+{{\Cal O}_{19}}="o19";%
(15,128)*+{{\Cal O}_{20}}="o20";%
(30,128)*+{{\Cal O}_{21}}="o21";%
(22,136)*+{{\Cal O}_{22}}="o22";%
(22,144)*+{{\Cal O}_{23}}="o23";%
(22,152)*+{{\Cal O}_{24}}="o24";%
(-15,0)*{0};
(-15,8)*{9};
(-15,16)*{12};
(-15,24)*{14};
(-15,32)*{15};
(-15,40)*{16};
(-15,48)*{17};
(-15,56)*{18};
(-15,64)*{19};
(-15,72)*{20};
(-15,80)*{21};
(-15,88)*{22};
(-15,96)*{23};
(-15,104)*{24};
(-15,112)*{25};
(-15,120)*{26};
(-15,128)*{27};
(-15,136)*{28};
(-15,144)*{29};
(-15,152)*{30};
{\ar@{-} "o0"; "o1"};%
{\ar@{-} "o1"; "o2"};%
{\ar@{-} "o1"; "o3"};%
{\ar@{-} "o2"; "o6"};%
{\ar@{-} "o3"; "o4"};%
{\ar@{-} "o3"; "o5"};%
{\ar@{-} "o4"; "o8"};%
{\ar@{-} "o5"; "o8"};%
{\ar@{-} "o6"; "o7"};%
{\ar@{-} "o7"; "o10"};%
{\ar@{-} "o7"; "o9"};%
{\ar@{-} "o8"; "o10"};%
{\ar@{-} "o9"; "o11"};%
{\ar@{-} "o9"; "o13"};%
{\ar@{-} "o9"; "o15"};%
{\ar@{-} "o10"; "o12"};%
{\ar@{-} "o10"; "o13"};%
{\ar@{-} "o10"; "o15"};%
{\ar@{-} "o11"; "o16"};%
{\ar@{-} "o13"; "o14"};%
{\ar@{-} "o13"; "o16"};%
{\ar@{-} "o14"; "o18"};%
{\ar@{-} "o15"; "o21"};%
{\ar@{-} "o15"; "o18"};%
{\ar@{-} "o16"; "o17"};%
{\ar@{-} "o14"; "o17"};%
{\ar@{-} "o20"; "o18"};%
{\ar@{-} "o19"; "o17"};%
{\ar@{-} "o19"; "o12"};%
{\ar@{-} "o19"; "o15"};%
{\ar@{-} "o19"; "o20"};%
{\ar@{-} "o21"; "o19"};%
{\ar@{-} "o21"; "o22"};%
{\ar@{-} "o22"; "o23"};%
{\ar@{-} "o23"; "o24"};%
{\ar@{-} "o20"; "o22"};%
\endxy
$$
\endproclaim

\bigskip\bigskip

$\spadesuit$. The hyperdiscriminant orbit closure $\overline{\OOO_{23}}$.

This is the hypersurface given by the tensors with vanishing hyperdiscriminant.
Its desingularization is, as always (see \cite{KW11a} section 5), a $1$-jet bundle $\xi$.
The orbit closure is not normal
but its normalization has rational singularities.

The minimal resolution $\FF(23)_\bullet$ of the normalization has terms
$$0\rightarrow (3,1,1;2^5)\rightarrow
(1^3;2,1^4)\oplus (0^3;0^5).$$
The extra representation in $\FF(24)_\bullet$ is just ${\goth g}_3$.

\bigskip\bigskip

$\spadesuit$ The codimension $2$ orbit closure $\overline{\OOO_{22}}$.

This is the cusp component of $Sing(\overline{\OOO_{3}})$ of codimension 2 in $\CC^3\otimes\bigwedge^3\CC^5$.
The resolution of the coordinate ring $\CC[ {\overline{\OOO_{22}}}]$ is determinantal, with the terms
$$(4^3;5^4,4)\rightarrow (4,4,2;4^5)\rightarrow (0^3;0^5).$$

The numerator of the Hilbert polynomial is
$5t^{10}+10t^9+9t^8+8t^7+7t^6+6t^5+5t^4+4t^3+3t^2+2t+1$.

The matrix in question can be viewed as a map from $S_2 E$ to $\bigwedge^4 F^*$ where the rows of the matrix are the
$4\times 4$ Pfaffians of our skew-symmetrc $5\times 5$ matrix, written as quadratic polynomials in three variables.

The orbit closure is not normal, as $\overline{\OOO_{21}}$ is contained in the singular locus
of $\overline{\OOO_{22}}$, as the rank of our matrix at a representative from $\OOO_{21}$ drops to three.

The desingularization of the normalization is given by a bundle $\eta$ whose maximal weights are $(1,0,0;0,0,0,1,1)$
and $(0,0,1;0,0,1,1,0)$. It lives on $Grass(2, E)\times Grass(1, F)$.
The resolution of the coordinate ring of the normalization is
$$(3,3,2;4,3^4)\rightarrow (3,2,2;3^4,2)\oplus (3,1,1;2^5)\rightarrow (1^3;2,1^4)\oplus (0^3;0^5).$$

The numerator of the Hilbert polynomial for the coordinate ring of normalization is
$15t^6+15t^5+15t^4+9t^3+3t^2+2t+1$.

\bigskip\bigskip

$\spadesuit$ The codimension $3$ orbit closure $\overline{\OOO_{21}}$.

The desingularization $Z(21)$ lives on $Grass(3, F)$. The bundle $\xi(21)$ for $Z(21)$ is $E\otimes\bigwedge^2{\Cal R}$, with $rank\ {\Cal R}=3$.
This orbit is normal, with rational singularities. The complex $\FF(21)_\bullet$ is as follows

$$0\rightarrow (3^3;4^3,3^2)\rightarrow (3,3,2;4,3^4)\rightarrow (3,1,1;2^5)\rightarrow (0^3;0^5).$$

Notice that this variety is the third secant of the orbit closure $\overline{\OOO_1}$, which is non-degenerate,  but there are no equations of degree $4$
vanishing on this secant. Its defining ideal is generated in degree $5$.

\bigskip\bigskip

For this representation, starting with this orbit closure we write the terms of {\it expected resolutions}. The exact shape is not certain
as there might be ''ghost" terms of pairs of cancelling representations. Federico Galetto is working on checking these resolutions
with Macaulay 2. We will update these results.

\bigskip\bigskip

$\spadesuit$ The codimension $3$ orbit closure $\overline{\Cal O}_{20}$.
The desingularization $Z$ corresponds to the bundle $\eta$
with weights
$$(1,0,0;1,1,0,0,0), (1,0,0;1,0,1,0,0), (1,0,0;1,0,0,1,0), (1,0,0;1,0,0,0,1),$$
$$(1,0,0;0,1,1,0,0),(1,0,0;0,1,0,1,0), (1,0,0;0,1,0,0,1), (1,0,0;0,0,1,1,0),$$
$$(0,1,0;1,1,0,0,0), (0,1,0;1,0,1,0,0), (0,1,0;1,0,0,1,0), (0,1,0;0,1,1,0,0),$$
$$(0,0,1;1,1,0,0,0), (0,0,1;1,0,1,0,0).$$

The corresponding complex $\FF(20)_\bullet$  is
$$0\rightarrow (5,3,2;4^5)\rightarrow (3,3,2;4,3^4)\oplus (4,2,2;4,3^4)\rightarrow$$
$$\rightarrow (3,1,1;2^5)\oplus (2^3;4,2^4)\oplus (3,2,2; 3^4,2)\rightarrow$$
$$\rightarrow (0^3;0^5)\oplus (1^3;2,1^4).$$

The orbit is not normal, but its normalization has rational singularities.

\bigskip\bigskip

$\spadesuit$ The codimension $4$ orbit closure $\overline{\Cal O}_{19}$.

The desingularization $Z$ corresponds to the bundle $\eta$
with weights
$$(1,0,0;1,1,0,0,0), (1,0,0;1,0,1,0,0), (1,0,0;1,0,0,1,0), (1,0,0;1,0,0,0,1),$$
$$(1,0,0;0,1,1,0,0),(1,0,0;0,1,0,1,0), (0,1,0;1,1,0,0,0), (0,1,0;1,0,1,0,0),$$
$$(0,1,0;1,0,0,1,0), (0,1,0;0,1,1,0,0), (0,0,1;1,1,0,0,0), (0,0,1;1,0,1,0,0).$$
$$(0,0,1;0,1,1,0,0).$$

The orbit closure is not normal but the normalization has rational singularities by Remark 1.1.

\proclaim{Conjecture}
The complex $\FF(19)_\bullet$ is
$$(5,4,4;6,5^4)\rightarrow(3^3;4^3,3^2)\oplus (5,3,2;4^5)\oplus (4,4,3;5^2,4^3)\oplus (4^3;5^4,4)\rightarrow$$
$$\rightarrow2*(3,3,2;4,3^4)\oplus (4,2,2;4,3^4)\oplus (3^3;4^4,2)\oplus (4,4,2;4^5)\rightarrow$$
$$\rightarrow 2*(3,1,1;2^5)\oplus (2^3;4,2^4)\oplus (3,2,2;3^4,2)\rightarrow(0^3;0^5)\oplus (1^3;2,1^4).$$
\endproclaim

The Euler characteristics of $\bigwedge^j\xi$ were calculated and they agree with the conjecture.

\bigskip\bigskip

$\spadesuit$ The codimension $6$ orbit closure $\overline{\Cal O}_{18}$.

The orbit closure is not normal, but normalization has rational singularities.
The desingularization $Z(18)$ lives on $Grass(2, E)\times Grass(4, F)$. The bundle $\xi$
is ${\Cal R}_E\otimes\bigwedge^2{\Cal R}_F$.
The terms of the complex $\FF(18)_\bullet$ are
$$(5,5,2;5^4,4)\rightarrow (5,4,2;5^2,4^3)\rightarrow (4,4,2;4^5)\oplus (5,2,2;4^3,3^2)\oplus (4,3,2;5,4,3^3)\rightarrow$$
$$(4,2,2;4^2,3^2,2)\oplus (4,2,2;4,3^4)\oplus (3,3,2;5,3^3,2)\oplus (4,1,1;3^2,2^3)\rightarrow$$
$$(3,2,2;4,3^3,1)\oplus (2^3;3^2,2^3)\oplus (4,1,0;2^5)\oplus (3,1,1;3,2^3,1)\rightarrow$$
$$(2^3;3^4,0)\oplus (2,1,1;2^4,0)\oplus (2,1,0;2,1^4)\rightarrow (1,1,0;1^4,0)\oplus (0^3;0^5).$$

\bigskip\bigskip

$\spadesuit$ The codimension $5$ orbit closure $\overline{\Cal O}_{17}$.

The desingularization $Z$ corresponds to the bundle $\eta$
with weights
$$(1,0,0;1,1,0,0,0), (1,0,0;1,0,1,0,0), (1,0,0;1,0,0,1,0), (1,0,0;1,0,0,0,1),$$
$$(1,0,0;0,1,1,0,0),(1,0,0;0,1,0,1,0), (0,1,0;1,1,0,0,0), (0,1,0;1,0,1,0,0),$$
$$(0,1,0;1,0,0,1,0), (0,1,0;0,1,1,0,0), (0,0,1;1,1,0,0,0), (0,0,1;1,0,1,0,0).$$

The orbit closure is normal with rational singularities by Remark 1.1.

\proclaim{Conjecture}
The complex $\FF(17)_\bullet$ has terms
$$(6,5,4;6^5)\rightarrow (5,4,4;6,5^4)\oplus (5,4,3;5^4,4)\rightarrow $$
$$(4,4,3;5^2,4^3)\oplus (5,3,2;4^5)\oplus (4,3,3;5,4^3,3)\oplus (3^3;4^3,3^2)\rightarrow$$
$$(3^3;4^4,2)\oplus (3^3;5,4,3^3)\oplus (4,2,2;4,3^4)\oplus (3,3,2;4,3^4)\rightarrow$$
$$(2^3;4,2^4)\oplus (3,1,1;2^5)\rightarrow (0^3;0^5).$$
\endproclaim

The Euler characteristics of $\bigwedge^j\xi$ were calculated and they agree with the conjecture.

\bigskip\bigskip

$\spadesuit$ The codimension $7$ orbit closure $\overline{\Cal O}_{16}$.

The desingularization $Z$ corresponds to the bundle $\eta$
with weights
$$(1,0,0;1,1,0,0,0), (1,0,0;1,0,1,0,0), (1,0,0;1,0,0,1,0), (1,0,0;1,0,0,0,1),$$
$$(1,0,0;0,1,1,0,0),(1,0,0;0,1,0,1,0), (1,0,0;0,1,0,0,1), (0,1,0;1,1,0,0,0),$$
$$(0,1,0;1,0,1,0,0), (0,1,0;1,0,0,1,0), (0,1,0;1,0,0,0,1),(0,1,0;0,1,1,0,0), $$
$$(0,1,0;0,1,0,1,0), (0,1,0;0,1,0,0,1), (0,0,1;1,1,0,0,0).$$

The orbit closure is normal with rational singularities by Remark 1.1.

\proclaim{Conjecture}
The complex $\FF(16)_\bullet$ has terms
$$(7,4,4;6^5)\rightarrow (6,4,3;6,5^4)\oplus (5,4,4;6,5^4)\rightarrow$$
$$(5,4,3;5^4,4)\oplus (4^3;5^4,4)\oplus (5,4,2;5^2,4^3)\oplus $$
$$(5,3,3;6,4^4)\oplus (4,4,3;6,4^4)\oplus (6,2,2;4^5)\rightarrow$$
$$(5,4,1;4^5)\oplus (4,4,2;5,4^3,3)\oplus (4,4,2;4^5)\oplus (4,3,3;5,4^3,3)\oplus$$
$$ (5,2,2;4^3,3^2)\oplus (4,3,2;5,4,3^3)\oplus (3^3;6,3^4)\rightarrow $$
$$(4,3,1;4,3^4)\oplus (4,2,2;4,4,3,3,2)\oplus (3,3,2;5,3,3,3,2)\oplus $$
$$(3,3,2;4,4,3,3,2)\oplus (3,2,2;3^4,2)\rightarrow$$
$$(3,2,2;4,3,3,3,1)\oplus (3,2,1;3,3,2,2,2)\oplus (2,2,2;3,3,3,2,1)\oplus (2,2,1;2^5)\rightarrow$$
$$(2^3;3^4,0)\oplus (2,1,1;2^3,1^2)\rightarrow (0^3;0^5).$$
\endproclaim

The Euler characteristics of $\bigwedge^j\xi$ were calculated and they agree with the conjecture.

\bigskip\bigskip

$\spadesuit$ The codimension $7$ orbit closure $\overline{\Cal O}_{15}$.

$$(1,0,0;1,1,0,0,0), (1,0,0;1,0,1,0,0), (1,0,0;1,0,0,1,0), (1,0,0;1,0,0,0,1),$$
$$(1,0,0;0,1,1,0,0), (0,1,0;1,1,0,0,0), (0,1,0;1,0,1,0,0), (0,1,0;1,0,0,1,0),$$
$$ (0,1,0;1,0,0,0,1),(0,1,0;0,1,1,0,0), (0,0,1;1,1,0,0,0), (0,0,1;1,0,1,0,0)$$
$$(0,0,1;0,1,1,0,0).$$

The orbit closure is not normal but the normalization has rational singularities by Remark 1.1.

The Euler characteristics of powers of $\xi$ are
$$\chi (\bigwedge^0\xi)= (0,0,0;0,0,0,0,0),$$
 $$\chi (\bigwedge^1\xi) = 0,$$
 $$\chi (\bigwedge^2\xi )= (1,1,0;1,1,1,1,0),$$
 $$\chi (\bigwedge^3\xi )= (2,1,0;2,1,1,1,1),$$
 $$\chi (\bigwedge^4\xi )= -(2,1,1;2,2,2,2,0)-(2,1,1;2,2,2,1,1),$$
 $$\chi (\bigwedge^5\xi )= -(4,1,0;2,2,2,2,2)-(3,1,1;3,2,2,2,1)-$$
 $$-(2,2,1;3,2,2,2,1)-(2,2,1;2,2,2,2,2),$$
 $$\chi (\bigwedge^6\xi)= -(4,1,1;3,3,2,2,2)-(2,2,2;3,3,3,3,0),$$
 $$\chi (\bigwedge^7\xi)= -(3,3,1;3,3,3,3,2)-(3,2,2;4,3,3,3,1),$$
$$\chi (\bigwedge^8\xi)= -(4,3,1;4,3,3,3,3)-(4,2,2;4,4,3,3,2)-(4,2,2;4,3,3,3,3)-$$
$$-2*(3,3,2;5,3,3,3,2)-(3,3,2;4,4,3,3,2)-(3,3,2;4,3,3,3,3),$$
$$\chi (\bigwedge^9\xi)= -(5,2,2;4,4,4,3,3)-2*(4,3,2;5,4,3,3,3)-(4,3,2;4,4,4,3,3)-$$
$$-(3,3,3;6,3,3,3,3)-(3,3,3;5,4,4,3,2)-(3,3,3;5,4,3,3,3),$$
$$\chi (\bigwedge^{10}\xi)= (4,4,2;4,4,4,4,4)-(4,3,3;5,5,4,3,3),$$
$$\chi (\bigwedge^{11}\xi)= (5,4,2;5,5,4,4,4)+(5,3,3;6,4,4,4,4),$$
$$\chi (\bigwedge^{12}\xi)= (5,5,2;5,5,5,5,4)-(4,4,4;6,6,4,4,4)-(4,4,4;5,5,5,5,4),$$
$$\chi (\bigwedge^{13}\xi)= -(5,4,4;6,6,5,5,4),$$
$$\chi (\bigwedge^{14}\xi)= -(5,5,4;6,6,6,6,4).$$

\bigskip\bigskip

$\spadesuit$. The codimension $7$ orbit closure $\overline{\Cal O}_{14}$.

$$(1,0,0;1,1,0,0,0), (1,0,0;1,0,1,0,0), (1,0,0;1,0,0,1,0), (1,0,0;1,0,0,0,1),$$
$$(1,0,0;0,1,1,0,0), (1,0,0;0,1,0,1,0), (1,0,0;0,1,0,0,1), (0,1,0;1,1,0,0,0),$$
$$(0,1,0;1,0,1,0,0), (0,1,0;0,1,1,0,0), (0,0,1;1,1,0,0,0), (0,0,1;1,0,1,0,0)$$
$$(0,0,1;0,1,1,0,0).$$

The orbit closure is not normal but the normalization has rational singularities by Remark 1.1.

The Euler characteristics of powers of $\xi$ are
$$\chi (\bigwedge^0\xi)= (0,0,0;0,0,0,0,0),$$
 $$\chi (\bigwedge^1\xi) = 0,$$
 $$\chi (\bigwedge^2\xi )= (1,1,0;1,1,1,1,0),$$
 $$\chi (\bigwedge^3\xi )= (2,1,0;2,1,1,1,1)+(1,1,1;2,1,1,1,1),$$
 $$\chi (\bigwedge^4\xi )= -(2,1,1;2,2,2,2,0),$$
 $$\chi (\bigwedge^5\xi )= -(4,1,0;2,2,2,2,2)-(3,1,1;3,2,2,2,1),$$
 $$\chi (\bigwedge^6\xi)= -(4,1,1;3,3,2,2,2)+(3,2,1;3,3,2,2,2)-(2,2,2;3,3,3,3,0)$$
 $$+(2,2,2;3,3,3,2,1)+(2,2,2;3,3,2,2,2),$$
 $$\chi (\bigwedge^7\xi)= (4,2,1;3,3,3,3,2)-(3,2,2;4,3,3,3,1)+$$
 $$+(3,2,2;4,3,3,2,2)+(3,2,2;3,3,3,3,2),$$
$$\chi (\bigwedge^8\xi)= -(4,2,2;4,4,3,3,2)-(3,3,2;5,3,3,3,2),$$
$$\chi (\bigwedge^9\xi)= -2*(5,2,2;4,4,4,3,3)-(4,3,2;5,4,3,3,3)-$$
$$-(3,3,3;6,3,3,3,3)-(3,3,3;4,4,4,3,3),$$
$$\chi (\bigwedge^{10}\xi)= -(6,2,2;4,4,4,4,4)+(4,4,2;4,4,4,4,4)+$$
$$+(4,3,3;5,4,4,4,3)+(4,3,3;4,4,4,4,4),$$
$$\chi (\bigwedge^{11}\xi)= (5,4,2;5,5,4,4,4)+(5,3,3;6,4,4,4,4)+(4,4,3;5,5,4,4,4),$$
$$\chi (\bigwedge^{12}\xi)= -(6,3,3;5,5,5,5,4)+(5,5,2;5,5,5,5,4),$$
$$\chi (\bigwedge^{13}\xi)= -(6,4,3;6,5,5,5,5),$$
$$\chi (\bigwedge^{14}\xi)= 0,$$
$$\chi (\bigwedge^{15}\xi)= (6,6,3;6,6,6,6,6).$$

\bigskip\bigskip

$\spadesuit$ The codimension $8$ orbit closure $\overline{\Cal O}_{13}$.

The desingularization $Z$ corresponds to the bundle $\eta$
with weights
$$(1,0,0;1,1,0,0,0), (1,0,0;1,0,1,0,0), (1,0,0;1,0,0,1,0), (1,0,0;1,0,0,0,1),$$
$$(1,0,0;0,1,1,0,0), (0,1,0;1,1,0,0,0), (0,1,0;1,0,1,0,0), (0,1,0;1,0,0,1,0),$$
$$ (0,0,1;1,1,0,0,0).$$

The orbit closure is not normal but the normalization has rational singularities by Remark 1.1.

\bigskip\bigskip

The Euler characteristics of powers of $\xi$ are
$$\chi (\bigwedge^0\xi)= (0,0,0;0,0,0,0,0),$$
 $$\chi (\bigwedge^1\xi) = 0,$$
 $$\chi (\bigwedge^2\xi )= (1,1,0;1,1,1,1,0),$$
 $$\chi (\bigwedge^3\xi )= (2,1,0;2,1,1,1,1)+(1,1,1;2,1,1,1,1),$$
 $$\chi (\bigwedge^4\xi )= -(2,1,1;2,2,2,2,0)-(2,1,1;2,2,2,1,1),$$
 $$\chi (\bigwedge^5\xi )= -(4,1,0;2,2,2,2,2)-(3,1,1;3,2,2,2,1)-(3,1,1;2,2,2,2,2)-$$
 $$-(2,2,1;3,2,2,2,1)-(2,2,1;2,2,2,2,2),$$
  $$\chi (\bigwedge^6\xi)= -(4,1,1;3,3,2,2,2)+(3,2,1;3,3,2,2,2)-$$
 $$-(2,2,2;3,3,3,3,0)+(2,2,2;3,3,3,2,1),$$
 $$\chi (\bigwedge^7\xi)= (4,2,1;3,3,3,3,2)-(3,2,2;4,3,3,3,1)+$$
 $$+(3,2,2;4,3,3,2,2)+(3,2,2;3,3,3,3,2),$$
$$\chi (\bigwedge^8\xi)= -(4,3,1;4,3,3,3,3)-(4,2,2;4,4,3,3,2)-2*(3,3,2;5,3,3,3,2)-$$
$$-(3,3,2;4,4,3,3,2)-(3,2,2;4,3,3,3,3),$$
$$\chi (\bigwedge^9\xi)= -2*(5,2,2;4,4,4,3,3)-2*(4,3,2;5,4,3,3,3)-(4,3,2;4,4,4,3,3)-$$
$$-2*(3,3,3;6,3,3,3,3)-(3,3,3;5,4,4,3,2)-(3,3,3;5,4,3,3,3),$$
$$\chi (\bigwedge^{10}\xi)= -(6,2,2;4,4,4,4,4)+(5,4,1;4,4,4,4,4)+(4,4,2;5,4,4,4,3)+$$
$$+(4,4,2;4,4,4,4,4)-(4,3,3;5,5,4,3,3)+$$
$$+(4,3,3;5,4,4,4,3)+(4,3,3;4,4,4,4,4),$$
$$\chi (\bigwedge^{11}\xi)= 2*(5,4,2;5,5,4,4,4)+2*(5,3,3;6,4,4,4,4)+$$
$$+(4,4,3;6,4,4,4,4)+(4,4,3;5,5,4,4,4),$$
$$\chi (\bigwedge^{12}\xi)= -(6,3,3;5,5,5,5,4)+(5,5,2;5,5,5,5,4)-(5,4,3;5,5,5,5,4)-$$
$$-(4,4,4;6,6,4,4,4)-(4,4,4;6,5,5,5,3)-(4,4,4;5,5,5,5,4),$$
$$\chi (\bigwedge^{13}\xi)= -2*(6,4,3;6,5,5,5,5)-2*(5,4,4;6,6,5,5,4)-(5,4,;6,5,5,5,5),$$
$$\chi (\bigwedge^{14}\xi)= -(6,4,4;6,6,6,5,5)-(5,5,4;6,6,6,6,4),$$
$$\chi (\bigwedge^{15}\xi)= (6,6,3;6,6,6,6,6)-(5,5,5;6,6,6,6,6),$$
$$\chi (\bigwedge^{16}\xi)= -(6,5,5;7,7,6,6,6),$$
$$\chi (\bigwedge^{17}\xi)= -(6,6,5;7,7,7,7,6),$$

\bigskip\bigskip

\head \S 7. The case $(E_7, \alpha_6)$. \endhead

 $X=E\otimes V(\omega_4 ,D_5 )$ where $E=\CC^2$,
$V(\omega_4 ,D_5 )$ is the half-spinor representation and ${ G}=SL(E)\times Spin
(10)\times \CC^*$.

The graded Lie algebra of type $E_7$ is
$${\goth g}(E_7)= {\goth g}_{-2}\oplus {\goth g}_{-1}\oplus {\goth g}_0\oplus {\goth g}_1\oplus {\goth g}_2$$
with $G_0= SL(2)\times Spin(10)\times\CC^*$, ${\goth g}_0= {\goth sl}(2)\oplus{\goth so}(10)\oplus \CC$,
${\goth g}_1=\CC^2\otimes V(\omega_4 ,D_5)$, ${\goth g}_2=\bigwedge^2\CC^2\otimes \CC^{10}$.

We label the weight vectors in $V(\omega_4;D_5 )$ by $[I]$ where $I$ is the subset of even cardinality of $\lbrace 1,2,3,4,5\rbrace$ where the sign of the component is negative. Thus the weight vectors in ${\goth g}_1$ are labelled by the pairs $[a;I]$ where $a\in\lbrace 1,2\rbrace$.

The invariant scalar product $(,)$ on $\goth g$ restricted to ${\goth g}_1$ is given by the formula
$$([a;I], [b;J]) = 1+\# (\lbrace a\rbrace\cap\lbrace b\rbrace )-{1\over 2}\# (\lbrace I\setminus J\rbrace\cup\lbrace J\setminus I\rbrace) .$$
Possible scalar products are $2,1,0,-1$.

$$\matrix number&{\goth s}&dim&representative\\
0&0&0&\\
1&A_1&12&[1;\emptyset ]\\
2&2A_1&17&[1;\emptyset ]+[1;1234]\\
3&2A_1&19&[1;\emptyset ]+[2;12]\\
4&3A_1&23&[1;\emptyset ]+[1;1234]+[2;12]\\
5&A_2&24&[1;\emptyset ]+[2;1234]\\
6&A_2+A_1&28&[1;\emptyset ]+[2;1234]+[1;1235]\\
7&A_2+2A_1&31&[1;\emptyset ]+[2;1234]+[1;1235]+[2;35] \\
8&2A_2&32&[1;\emptyset ]+[2;1234]+[1;1235]+[2;45] \\
\endmatrix$$

$$\matrix number&description&geometry\\
0&0&\\
1&h.weight\ vector&\\
2&E-degenerate&\\
3&V(\omega_4;D_5 )-degenerate&\\
4&&\tau (\overline\OOO_1)\\
5&&\sigma_2 (\overline\OOO_1)\\
6&\exists\ a\ member\ which\ is\ pure&\\
7&hyperdiscriminant& \\
8&general& \\
\endmatrix$$

\bigskip

The numerical data is as follows

$$\matrix number&degree&numerator\\
0&1&1\\
1&132&1+20t+60t^2+44t^3+7t^4\\
2&16&1+15t\\
3&408&1+13t+61t^2+129t^3+129t^4+61t^5+13t^6+t^7\\
4&584&1+9t+45t^2+133t^3+201t^4+145t^5+45t^6+5t^7\\
5&388&1+8t+36t^2+88t^3+122t^4+88t^5+36t^6+\\
&&+8t^7+t^8\\
6&60&1+4t+20t^2+28t^3+7t^4\\
7&4&1+t+t^2+t^3\\
8&1&1
\endmatrix$$

The singularities data is

  $$\matrix number&spherical&normal&C-M&R.S.&Gor\\
0&yes&yes&yes&yes&yes\\
1&yes&yes&yes&yes&no\\
2&yes&yes&yes&yes&no\\
3&yes&yes&yes&yes&yes\\
4&yes&yes&yes&yes&no\\
5&no&yes&yes&yes&yes\\
6&no&no&yes&no&no\\
n(6)&no&yes&yes&yes&no\\
7&no&yes&yes&yes&yes\\
8&no&yes&yes&yes&yes
\endmatrix$$

\proclaim{Remark} The degeneration order is
$$
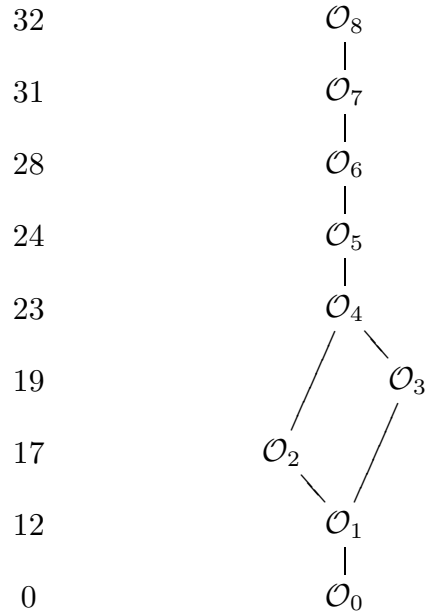
\xy
(20,0)*+{{\Cal O}_{0}}="o0";%
(20,8)*+{{\Cal O}_{1}}="o1";%
(13,16)*+{{\Cal O}_2}="o2";%
(27,24)*+{{\Cal O}_3}="o3";%
(20,32)*+{{\Cal O}_{4}}="o4";%
(20,40)*+{{\Cal O}_{5}}="o5";%
(20,48)*+{{\Cal O}_{6}}="o6";%
(20,56)*+{{\Cal O}_{7}}="o7";%
(20,64)*+{{\Cal O}_{8}}="o8";%
(-15,0)*{0};
(-15,8)*{12};
(-15,16)*{17};
(-15,24)*{19};
(-15,32)*{23};
(-15,40)*{24};
(-15,48)*{28};
(-15,56)*{31};
(-15,64)*{32};
{\ar@{-} "o0"; "o1"};%
{\ar@{-} "o2"; "o1"};%
{\ar@{-} "o3"; "o1"};%
{\ar@{-} "o2"; "o4"};%
{\ar@{-} "o3"; "o4"};%
{\ar@{-} "o4"; "o5"};%
{\ar@{-} "o5"; "o6"};%
{\ar@{-} "o6"; "o7"};%
{\ar@{-} "o7"; "o8"};%
\endxy
$$
\endproclaim

$\spadesuit$ The hyperdiscriminant orbit closure $\overline{\OOO_{7}}$.

This is the hypersurface given by the tensors with vanishing hyperdiscriminant.
Its desingularization is, as always (see \cite{KW11a} section 5), a $1$-jet bundle $\xi$.
The orbit closure is normal.

The minimal resolution $\FF(8)_\bullet$ of the normalization has terms
$$0\rightarrow (2,2;0^5)\rightarrow
 (0,0;0^5).$$

\bigskip\bigskip

$\spadesuit$ The codimension $4$ orbit closure $\overline{\OOO_{6}}$.

This is the orbit closure of pencils of spinors containing a pure spinor.
Its resolution is obtained from the resolution of the variety of pure spinors in $V(\omega_4 ,D_5)$ in the same way
as analogous cases for the determinantal  varieties.
The bundle $\eta$ has weights

$$ (1,0;+,+,+,+,+), (1,0;+,+,+,-,-), (1,0;+,+,-,+,-), $$
$$(1,0;+,-,+,+,-),(1,0;+,+,-,-,+), (1,0;-,+,+,+,-),$$
$$(0,1;+,+,+,+,+), (0,1;+,+,+,-,-), (0,1;+,+,-,+,-)$$
$$(0,1;+,-,+,+,-).$$

The terms of the complex $\FF(7)_\bullet$ are
$$(7,1;0^5)\rightarrow (5,1;1,0^4)\rightarrow (4,1;\omega_4 )\rightarrow$$
$$\rightarrow (2,1;\omega_5)\rightarrow (1,1;1,0^4)\oplus (0,0;0^5).$$

The extra representation in $\FF(7)_0$ is ${\goth g}_2$ so, as always, we can identify the normalization as an orbit closure in ${\goth g}_1\oplus{\goth g}_2$.

\bigskip\bigskip

$\spadesuit$  The codimension $8$ orbit closure $\overline{\OOO_{5}}$.

The secant $\sigma(\overline{\OOO_1})$ of the orbit closure of the highest weight vector. This orbit closure has rational singularities, is Gorenstein, of codimension 8. The defining ideal is generated by the representation $S_{2,1}E\otimes V(\omega_4 )$ in degree 3. The bundle $\xi =E\otimes\xi'$ where the bundle $\xi'$ has weights with ${-1}\over 2$ on the first coordinate. The desingularization $Z(6)$ lives over the isotropic Grassmannian $IGrass(1, \CC^{10} )$.

The orbit closure is normal with rational singularities by Remark 1.1.

\proclaim{Conjecture}
The terms in the finite free resolution of the coordinate ring are:
$$0\rightarrow S_{8,8}E\otimes A(-16)\rightarrow S_{7,6}E\otimes V_{\omega_4}\otimes A(-13)\rightarrow $$
$$\rightarrow  S_{7,5}E \otimes A(-12)\oplus S_{6,6}E\otimes V_{\omega_2}\otimes A(-12)\oplus S_{7,4}E\otimes V_{\omega_5}\otimes A(-11) \rightarrow $$
$$\rightarrow S_{7,3}E\otimes V_{\omega_1}\otimes A(-10)\oplus S_{6,4}E\otimes V_{\omega_1}\otimes A(-10)\oplus S_{5,4}E\otimes V_{\omega_1+\omega_5}\otimes A(-9) \rightarrow$$
$$\rightarrow S_{7,1}E\otimes A(-8)\oplus S_{6,2}E\otimes A(-8)\oplus S_{5,3}E\otimes V_{\omega_2}\otimes A(-8)\oplus S_{5,3}E\otimes V_{2\omega_1}\otimes A(-8)$$
$$\oplus S_{4,4}E\otimes A(-8)\oplus S_{4,4}E\otimes V_{\omega_4+\omega_5}\otimes A(-8)\oplus S_{4,4}E\otimes V_{2\omega_1}\otimes A(-8)\rightarrow$$
$$\rightarrow S_{4,3}E\otimes V_{\omega_1+\omega_4}\otimes A(-7)\oplus S_{5,1}E\otimes V_{\omega_1}\otimes A(-6)\oplus S_{4,2}E\otimes V_{\omega_1}\otimes A(-6)\rightarrow$$
$$\rightarrow S_{4,1}E\otimes V_{\omega_4}\otimes A(-5)\oplus  S_{3,1}E \otimes A(-4)\oplus S_{2,2}E\otimes V_{\omega_2}\otimes A(-4)\rightarrow $$
$$\rightarrow S_{2,1}E\otimes V_{\omega_5}\otimes A(-3)\rightarrow A$$
\endproclaim

The Euler characteristics of $\bigwedge^j\xi$ were calculated and they agree with the conjecture.

\bigskip\bigskip

$\spadesuit$ The codimension $9$ orbit closure $\overline{\OOO_{4}}$.

This orbit closure is the tangential variety of the highest weight orbit closure $\bar{\OOO_1}$.

The orbit closure is  normal with rational singularities by Remark 1.1.

The Euler characteristics of the exterior powers of $\xi$ are as follows.

$$\chi (\bigwedge^0\xi )=S_{0,0}E,$$
$$\chi (\bigwedge^1\xi )=0,$$
$$\chi (\bigwedge^2\xi )=0,$$
$$\chi (\bigwedge^3\xi )=S_{2,1}E\otimes V_{\omega_5},$$
$$\chi (\bigwedge^4\xi )=S_{2,2}E\otimes V_{\omega_2}-S_{2,2}E\otimes V_{2\omega_1}+S_{3,1}E,$$
$$\chi (\bigwedge^5\xi )=-S_{3,2}E\otimes V_{\omega_1+\omega_5}-S_{4,1}E\otimes V_{\omega_4},$$
$$\chi (\bigwedge^6\xi )=-S_{4,2}E\otimes V_{\omega_3}- S_{3,3}E\otimes V_{\omega_1}-S_{4,2}E\otimes V_{\omega_1}-$$
$$-S_{5,1}E\otimes V_{\omega_1}-S_{3,3}E\otimes V_{2\omega_5},$$
$$\chi (\bigwedge^7\xi )=-S_{4,3}E\otimes V_{\omega_5}-S_{4,3}E\otimes V_{\omega_1+\omega_4},$$
$$\chi (\bigwedge^8\xi )=S_{5,3}E\otimes V_{2\omega_1}+S_{5,3}E\otimes V_{\omega_4+\omega_5}+S_{5,3}E\otimes V_{\omega_2}+$$
$$+ S_{6,2}E\otimes V_{\omega_2}+S_{4,4}E\otimes V_{\omega_4+\omega_5}+S_{6,2}E+$$
$$+ S_{4,4}E\otimes V_{2\omega_1}+S_{7,1}E,$$
$$\chi (\bigwedge^9\xi )=S_{6,3}E\otimes V_{\omega_1+\omega_5}+S_{5,4}E\otimes V_{\omega_1+\omega_5}+$$
$$+S_{6,3}E\otimes V_{\omega_4}+S_{7,2}E\otimes V_{\omega_4},$$
$$\chi (\bigwedge^{10}\xi )= -S_{6,4}E\otimes V_{2\omega_4}-S_{6,4}E\otimes V_{\omega_1}-S_{5,5}E\otimes V_{\omega_3},$$
$$\chi (\bigwedge^{11}\xi )=-S_{7,4}E\otimes V_{\omega_1+\omega_4}-S_{6,5}E\otimes V_{\omega_1+\omega_4}-S_{8,3}E\otimes V_{\omega_5}-$$
$$-S_{7,4}E\otimes V_{\omega_5}-S_{6,5}E\otimes V_{\omega_5},$$
$$\chi (\bigwedge^{12}\xi )=-S_{7,5}E\otimes V_{2\omega_1}-S_{8,4}E\otimes V_{\omega_2}-S_{9,3}E,$$
$$\chi (\bigwedge^{13}\xi )= S_{8,5}E\otimes V_{\omega_4}+ S_{7,6}E\otimes V_{\omega_4},$$
$$\chi (\bigwedge^{14}\xi )= S_{9,5}E\otimes V_{\omega_1}+S_{8,6}E\otimes V_{\omega_1},$$
$$\chi (\bigwedge^{15}\xi )= 0,$$
$$\chi (\bigwedge^{16}\xi )=-S_{10,6}E.$$

The orbit closure is normal. The defining ideal is generated by cubics (equations of the secant $\overline{\OOO_5}$) and quartics.

\bigskip\bigskip

\head \S 8. The case $(E_7, \alpha_7)$. \endhead

$X= V(\omega_6 ,E_6 )$, the sixth fundamental
representation for the group $G(E_6)$, with $G_0=G(E_6)\times\CC^*$ where $G(E_6)$ is a simply connected group of type $E_6$.

The graded Lie algebra of type $E_7$ is
$${\goth g}(E_7)=  {\goth g}_{-1}\oplus {\goth g}_0\oplus {\goth g}_1$$
with ${\goth g}_0=  {\goth g}(E_6)\oplus\CC$,
${\goth g}_1= V(\omega_6 , E_6)$.

The weight vectors of $V(\omega_6 ,E_6)$ are parametrized by the roots of $E_7$ whose coefficient
of $\alpha_7$ equals $1$. There are 27 such roots; we index them by labelled Dynkin diagram
(with the coefficient of $\alpha_7$  being $1$).
The roots are
$$\matrix &&0&&&\\ 0&0&0&0&0&1\endmatrix$$
$$\matrix &&0&&&\\ 0&0&0&0&1&1\endmatrix$$
$$\matrix &&0&&&\\ 0&0&0&1&1&1\endmatrix$$
$$\matrix &&0&&&\\ 0&0&1&1&1&1\endmatrix$$
$$\matrix &&0&&&\\ 0&1&1&1&1&1\endmatrix ,\ \  \matrix &&1&&&\\ 0&0&1&1&1&1\endmatrix$$
$$\matrix &&0&&&\\ 1&1&1&1&1&1\endmatrix ,\ \  \matrix &&1&&&\\ 0&1&1&1&1&1\endmatrix$$
$$\matrix &&1&&&\\ 1&1&1&1&1&1\endmatrix ,\ \  \matrix &&1&&&\\ 0&1&2&1&1&1\endmatrix$$
$$\matrix &&1&&&\\ 1&1&2&1&1&1\endmatrix ,\ \  \matrix &&1&&&\\ 0&1&2&2&1&1\endmatrix$$
$$\matrix &&1&&&\\ 1&2&2&1&1&1\endmatrix ,\  \ \matrix &&1&&&\\ 1&1&2&2&1&1\endmatrix ,\ \ \matrix &&1&&&\\ 0&1&2&2&2&1\endmatrix$$
$$\matrix &&1&&&\\ 1&2&2&2&1&1\endmatrix ,\ \ \matrix &&1&&&\\ 1&1&2&2&2&1\endmatrix $$
$$\matrix &&1&&&\\ 1&2&3&2&1&1\endmatrix ,\  \ \matrix &&1&&&\\ 1&2&2&2&2&1\endmatrix $$
$$\matrix &&2&&&\\ 1&2&3&2&1&1\endmatrix ,\ \  \matrix &&1&&&\\ 1&2&3&2&2&1\endmatrix $$
$$\matrix &&2&&&\\ 1&2&3&2&2&1\endmatrix ,\ \ \matrix &&1&&&\\ 1&2&3&3&2&1\endmatrix $$
$$\matrix &&2&&&\\ 1&2&3&3&2&1\endmatrix $$
$$\matrix &&2&&&\\ 1&2&4&3&2&1\endmatrix $$
$$\matrix &&2&&&\\ 1&3&4&3&2&1\endmatrix $$
$$\matrix &&2&&&\\ 2&3&4&3&2&1\endmatrix $$

There is one invariant $\Delta$ of degree $3$ and the orbits are as follows

$$\matrix number&{\goth s}&representative&dim&spherical&normal&C-M&R.S.&Gor\\
0&0&&0&yes&yes&yes&yes&yes\\
1&A_1&&17&yes&yes&yes&yes&yes\\
2&2A_1&&26&yes&yes&yes&yes&yes\\
3&3A_1&&27&yes&yes&yes&yes&yes\\
\endmatrix$$

The numerical data are as follows

$$\matrix number&degree&numerator\\
0&1&1\\
1&78&1+10t+28t^2+28t^3+10t^4+t^5\\
2&3&1+t+t^2\\
3&1&1
\endmatrix$$

\proclaim{Remark} The degeneration order is linear.
\endproclaim

In order to exhibit representatives we just need to show three ortogonal weight vectors.
They are
$$\matrix &&0&&&\\ 0&0&0&0&0&1\endmatrix ,\  \matrix &&1&&&\\ 0&1&2&2&2&1\endmatrix ,\ \matrix &&2&&&\\ 2&3&4&3&2&1\endmatrix .$$
The invariant $\Delta$ of degree $3$ is just
$$\Delta =\sum_{\lbrace \gamma_1 ,\gamma_2 ,\gamma_3\rbrace =3A_1}\pm x_{\gamma_1}x_{\gamma_2}x_{\gamma_3}$$

We denote $[a,b,c,d,e,f]$ the highest weight module for $E_6$ with the highest weight
$a\omega_1+b\omega_2+\ldots +f\omega_6$.

$\spadesuit$
The resolution of $\CC[\bar{\OOO_1}]$ is as follows
$$0\rightarrow [0,0,0,0,0,0](-15)\rightarrow [1,0,0,0,0,0](-13)\rightarrow [0,1,0,0,0,0](-12)\rightarrow$$
$$[0,0,0,0,1,0](-10)\rightarrow [1,0,0,0,0,1](-9)\rightarrow [2,0,0,0,0,0](-8)\oplus [0,0,0,0,0,2](-7)$$
$$\rightarrow[1,0,0,0,0,1](-6)\rightarrow [0,0,1,0,0,0](-5)\rightarrow [0,1,0,0,0,0](-3)\rightarrow$$
$$[0,0,0,0,0,1](-2)\rightarrow [0,0,0,0,0,0].$$
This fact was verified in \cite{G11} and in \cite{SW13}.

\bigskip\bigskip

\head \S 9. Conclusions. \endhead

Here are some general conclusions we checked type by type.

Let $X_n$ be the Dynkin diagram of type $E_7$. The orbit  $X_{disc.}$ has the following properties.

\proclaim {Proposition 9.1}
\item {a)} The variety $X_{disc.}$ is closed $G_0$-equivariant and irreducible and thus it is an orbit closure,
\item{b)} The variety $X_{disc.}$ is a hypersurface if and only if the ring of invariants \break $Sym({\goth g}_1^*)^{(G,G)}$ contains a non-constant invariant. In these cases we have
$$Sym({\goth g}_1^*)^{(G,G)}=K[\Delta ]$$
 and $X_{disc.}$ is a hypersurface given by vanishing of the invariant $\Delta$,
\item{c)} In the cases when $Sym({\goth g}_1^*)^{(G,G)}=K$ the variety $X_{disc.}$ has codimension bigger than one.
\endproclaim

\proclaim {Proposition 9.2}
Let $X_n$ be of type $E_7$.
\item{a)} The orbit closure is spherical if and only if the support algebra $\goth s$ has all simple components of type $A_1$,
\item{b)} If the orbit closure $\overline{\OOO_v}$ is not normal then its normalization is contained in the representation ${\goth g}_1\oplus {\goth g}_i$ for some $i>1$,
where ${\goth g}_i$ is the $i$-th graded component in the grading associate to the simple root $\alpha_k$.
\endproclaim



\Refs\widestnumber\key{ABW82}



\ref
\key ABW82
\by Akin, K., Buchsbaum, D., Weyman, J.
\paper Schur Functors and Schur Complexes
\jour Adv. in Math.
\vol 44
\yr 1982
\pages 207-278
\endref

\ref
\key BC76a
\by Bala, P., Carter, R.
\paper Classes of Unipotent Elements in Simple Algebraic Groups I,
\jour Proc. Camb. Phil. Soc.
\vol 79
\yr 1976
\pages 401-425
\endref

\ref
\key BC76b
\by Bala, P., Carter, R.
\paper Classes of Unipotent Elements in Simple Algebraic Groups II,
\jour Proc. Camb. Phil. Soc.
\vol 80
\yr 1976
\pages 1-18
\endref

\ref
\key BE73
\by Buchsbaum, D. A.; Eisenbud, D.,
\paper What makes a complex exact?
\jour J. Algebra
\vol 25
\yr1973
\pages  259--268
\endref


\ref
\key CM93 \by Collingwood, D., McGovern, W. \book Nilpotent Orbits in Semisimple Lie Algebras \eds
\publ Van Nostrand Reinhold \publaddr New York
\yr 1993 \finalinfo Van Nostrand Reinhold Mathematics Series
\endref

\ref
\key DK85
 \by J.~Dadok, V.~Kac,
 \paper Polar Representations,
 \jour
 J. of Algebra,
 \vol 92
 \yr 1985
 \pages 504-524
 \endref

\ref
\key dG11 \by de Graaf, W.
\book  SLA - a GAP package \eds
\publ available at http://www.science.unitn.it/~degraaf/ \publaddr
\yr 2011 \finalinfo \endref

\ref
\key EK94 \by Eisebud ,D., Koh J.
\paper Nets of alternating matrices and Linear Syzygy Conjecture
\jour Adv. in Math.
\vol 106
\yr 1994
\pages 1-35
\endref

\ref
\key FH91 \by Fulton, W.; Harris, J. \book Representation Theory \eds
\publ Springer-Verlag \publaddr New York Heidelberg Berlin
\yr 1991 \finalinfo Graduate Texts in Math. vol. 129
\endref

\ref
\key G11
\by Galetto, F.,
\paper Resolutions of coordinate rings of orbit closures for representations with finitely many orbits
\jour
\vol
\yr
\pages
\endref




\ref
\key K82
\by Kac, V.
\paper Some remarks on nilpotent orbits
\jour J. of Algebra
\vol 64
\yr 1982
\pages 190-213
\endref


\ref
\key KW11a
\by Kra\'skiewicz, W., Weyman, J.
\paper Geometry of orbit closures for the representations associated to gradings of Lie algebras of types $E_6$, $F_4$ and $G_2$.
\jour
\vol
\yr
\pages
\endref

\ref
\key LM01
\by Landsberg, J., Manivel, L.
\paper The projective geometry of Freudenthal magic square
\jour J. of Algebra
\vol 239(2)
\yr 2001
\pages 477-512
\endref

\ref
\key P01
\by Parfenov, P.G.,
\paper Orbits and their closures in $\CC^{k_1}\otimes\ldots\otimes \CC^{k_r}$
\jour Matematicheskii Sbornik
\vol 192 (1)
\yr 2001
\pages 89-111
\endref


\ref
\key PW86
\by Pragacz, P., Weyman, J.
\paper On the construction of resolutions of determinantal varieties; a survey
\book Lecture Notes in Mathematics
\vol 1220
\yr 1986
\pages 73-92
\endref

\ref
\key R10 \by Ribeiro, J.
\book  Roots and weights python package \eds
\publ  http://www.math.neu.edu/~weyman/ \publaddr
\yr 2010 \finalinfo \endref


\ref
\key SW13
\by Sam, S.V, Weyman, J.
\paper Littlewood-type complexes and analogues of determinantal varieties
\publ http://www.math.neu.edu/~weyman/ \publaddr
\yr 2013 \finalinfo \endref


\ref
\key T06 \by Timashev, D.A.
\book  Homogeneous spaces and equivariant embeddings
\eds to appear
\publ arxiv:math/0602228
\publaddr
\yr 2006 \finalinfo \endref


\ref
\key V75 \by Vinberg, E.B.
\paper Weyl group of a graded Lie algebra
\jour Izv. Akad. Nauk SSSR
\vol 40
\yr 1975
\pages 488-526
\endref

\ref
\key V87 \by Vinberg, E.B.
\paper Classification of homogeneous nilpotent elements
of a semisimple graded Lie algebra
\jour Selecta Mathematica Sovietica
\vol 6 no.1
\yr 1987
\endref



\ref
\key W03 \by Weyman, J.
\book Cohomology of vector bundles and syzygies
\publ Cambridge University Press
\publaddr Cambridge, UK
\eds \yr 2003 \pages \finalinfo
Cambridge Tracts in Mathematics, vol. 149
\endref

\ref
\key WZ96 \by Weyman, J., Zelevinsky, A.
\paper Singularities of hyperdeterminants
\jour Ann. Inst. Fourier (Grenoble)
\vol 46
\year 1996
\pages 591--644
\endref
\endRefs
\enddocument